%% file: main.tex
\documentclass[11pt]{article}

\usepackage[utf8]{inputenc} 
\usepackage[T1]{fontenc}    

\input{command.tex}

\usepackage{enumerate}
\usepackage{natbib}

\usepackage{xr}

\newtheorem{claim}{Claim}[section]
\newtheorem{lemma}[claim]{Lemma}
\newtheorem{assumption}{Assumption}[section]
\newtheorem{definition}{Definition}[section]
\newtheorem{theorem}{Theorem}[section]

\newtheorem{corollary}{Corollary}[section]

\author{
Yang Peng\thanks{Yau Mathematical Sciences Center, Tsinghua University; email: \texttt{yang-peng@mail.tsinghua.edu.cn}.} \and
Yuchen Xin\thanks{School of Mathematical Sciences, Peking University; email: \texttt{2301110087@pku.edu.cn}.} \and
Zhihua Zhang\thanks{School of Mathematical Sciences, Peking University; email: \texttt{zhzhang@math.pku.edu.cn}.}
}

\title{Matrix Probability Inequalities under the Spectral Norm for Martingales and Ergodic Markov Chains with Applications in Statistical Learning}

\begin{document}
\maketitle

\begin{abstract}
In this paper, we study moment and concentration inequalities for the spectral norm of sums of dependent random matrices. 
We establish novel Rosenthal-Burkholder inequalities for the spectral norm of discrete-time matrix local martingales, Burkholder-Davis-Gundy inequality for the spectral norm of  continuous matrix local martingales, as well as matrix Rosenthal, Hoeffding, and Bernstein inequalities under the spectral norm for ergodic Markov chains. 
Compared with previous work on matrix concentration inequalities for Markov chains, which assume a non-zero absolute $L^2$-spectral gap or the stronger $\psi$-mixing condition, our results assume geometric ergodicity, a condition commonly used in statistical applications.
Furthermore, for uniformly geometrically ergodic chains, our bounds have leading terms that match the Markov chain central limit theorem.
We also give dimension-free versions of the martingale inequalities, which depend on the effective rank instead of the ambient dimension $d$ and extend to linear operators in infinite-dimensional Hilbert spaces. 
Our results have extensive applications in statistics and machine learning; in particular, we obtain improved bounds in covariance estimation and principal component analysis on Markovian data.
\end{abstract}

\section{Introduction}\label{Section:intro}
\input{intro}

\section{Main Results}\label{Section:main_results}
\input{main_results_martingale}
\input{main_results_markov}

\input{main_results_exteinsion_V_ergo}

\section{Applications}\label{Section:application}
\input{application}

\section{Proofs}\label{Section:proof_outline}
\input{proof_outline_burkholder}
\input{proof_outline_markov}

\section{Concluding Remarks}\label{Section:conclusion}
\input{conclusion}

\begin{appendix}
\section{Omitted Proof for the Inequalities for Matrix Martingales}\label{Section:proof_burkholder}
\input{proof_burkholder}

\section{Omitted Proof for the Inequalities for Markov Chains}\label{Section:proof_markov}
\input{proof_markov}

\end{appendix}

\bibliographystyle{abbrvnat}
\bibliography{ref}       

\end{document}

%% file: command.tex
\usepackage{geometry}
\usepackage{setspace}
\usepackage{amsmath, amssymb, amsfonts, bm, mathtools, mathrsfs}
\usepackage{amsthm}
\usepackage[dvipsnames]{xcolor}
\definecolor{darkblue}{rgb}{0,0,.5}
\usepackage{graphicx}
\usepackage{subfigure}
\usepackage[numbers]{natbib}
\usepackage[colorlinks=true,allcolors=darkblue]{hyperref}       
\usepackage{url}            
\usepackage{booktabs}       
\usepackage{amsfonts}       
\usepackage{nicefrac}       
\usepackage{microtype}      

\allowdisplaybreaks

\usepackage{float}
\usepackage{multirow}
\usepackage{footnote}
\usepackage{dsfont}
\usepackage{mathabx}

\usepackage{algorithm}
\usepackage{algorithmic}
\usepackage{nicefrac}
\usepackage{tikz}
\usepackage{overpic}

\usepackage{dsfont}
\usepackage{hyperref}
\usepackage[capitalize]{cleveref}
\usepackage{crossreftools}

\makeatletter
\newcommand*{\rom}[1]{\expandafter\@slowromancap\romannumeral #1@}
\makeatother

\newcommand{\ind}{\mathds{1}}

\newcommand{\brc}[1]{\left\{{#1}\right\}}
\newcommand{\prn}[1]{\left({#1}\right)} 
\newcommand{\brk}[1]{\left[{#1}\right]} 
\newcommand{\norm}[1]{\left\|{#1}\right\|} 
\newcommand{\abs}[1]{\left|{#1}\right|} 
\newcommand{\qv}[1]{\left\langle{#1}\right\rangle}

\newcommand{\wtilde}[1]{\widetilde{#1}}

\def\rQ{{\mathrm{Q}}}

\def\gF{{\mathcal{F}}}

\def\gP{{\mathcal{P}}}

\def\gX{{\mathcal{X}}}

\def\gZ{{\mathcal{Z}}}

\def\sF{{\mathscr{F}}}
\def\sG{{\mathscr{G}}}
\def\sH{{\mathscr{H}}}
\def\sZ{{\mathscr{Z}}}

\def\bB{{\bm{B}}}

\def\bD{{\bm{D}}}
\def\bh{{\bm{h}}}
\def\bH{{\bm{H}}}
\def\bN{{\bm{N}}}
\def\bM{{\bm{M}}}

\def\bF{{\bm{F}}}
\def\bW{{\bm{W}}}

\def\bI{{\bm{I}}}

\def\bG{{\bm{G}}}

\def\bK{{\bm{K}}}

\def\bPhi{{\bm{\Phi}}}

\def\bSigma{{\bm{\Sigma}}}

\def\bu{{\bm{u}}}
\def\bh{{\bm{h}}}
\def\bv{{\bm{v}}}

\def\bV{{\bm{V}}}
\def\bY{{\bm{Y}}}
\def\bX{{\bm{X}}}
\def\bR{{\bm{R}}}

\def\bD{{\bm{D}}}
\def\bS{{\bm{S}}}

\def\bUp{{\bm{\Upsilon}}}

\def\tmix{t_{\operatorname{mix}}}

\def\CB{{\mathbb C}}
\def\RB{{\mathbb R}}
\def\EB{{\mathbb E}}

\def\PB{{\mathbb P}}

\def\NB{{\mathbb N}}
\def\HB{{\mathbb H}}

\newcommand{\TV}{\textup{TV}}

\makeatletter
\long\def\@makecaption#1#2{
  \vskip 0.8ex
  \setbox\@tempboxa\hbox{\small {\bf #1:} #2}
  \parindent 1.5em  
  \dimen0=\hsize
  \advance\dimen0 by -3em
  \ifdim \wd\@tempboxa >\dimen0
  \hbox to \hsize{
    \parindent 0em
    \hfil 
    \parbox{\dimen0}{\def\baselinestretch{0.96}\small
      {\bf #1.} #2
    } 
    \hfil}
  \else \hbox to \hsize{\hfil \box\@tempboxa \hfil}
  \fi
}
\makeatother

\newcommand{\tr}{\operatorname{Tr}}

%% file: intro.tex
Moment and concentration inequalities for the spectral norm of sums of random matrices \citep{wigderson2005randomness,tropp2015introduction} are of great significance for statistics, machine learning, theoretical computer science,  signal processing, etc.
Common examples include covariance matrix estimation \citep{chen2012masked}, principal component analysis (PCA) \citep{kumar2023streaming}, and reinforcement learning  \citep{duan2022policy,wu2025uncertainty,mou2025statistical}.
Although matrix probability inequalities  have been studied extensively in the literature \citep{oliveira2009concentration,tropp2011freedmans, tropp2012user, mackey2014matrix,tropp2015introduction,paulin2016efron,banna2016bernstein, minsker2017some,tropp2018second,garg2018matrix, bacry2018concentration,klochkov2020uniform,qiu2020matrix,han2020moment,huang2021poincare,cai2022non,bandeira2023matrix,zhivotovskiy2024dimension,neeman2024concentration,kroshnin2024bernstein,brailovskaya2024universality,jirak2024concentration,wang2024sharp,wu2025uncertainty,van2025matrix}, the Rosenthal-Burkholder-type inequalities for the spectral norm of matrix martingales and matrix functions of Markov chains remain less explored. 
While \citet{dirksen2023tuning} recently established Rosenthal-Burkholder-type inequalities for the spectral norm of matrix martingales, the dependence of their constants on the moment order is suboptimal. Moreover, to the best of our knowledge, the Rosenthal-type inequalities for the spectral norm of matrix functions of Markov chains remain absent.
 
We bridge the gaps in the current paper.
Specifically, we first give a Rosenthal-Burkholder inequality \citep{rosenthal1970subspaces,burkholder1973distribution,hitczenko1990best} for the spectral norm of discrete-time matrix martingales in Theorem~\ref{thm:martingale_rosenthal}.
That is, for an $(\sF_k)_{k\in\NB}$-adapted $d$-dimensional self-adjoint matrix martingale $(\bM_k)_{k\in\NB}$, letting martingale differences $\bX_k=\bM_k-\bM_{k-1}$ and predictable quadratic variation $\qv{\bM}_{\infty}=\sum_{k=0}^\infty \EB[\bX_k^2|\sF_{k-1}]$, one has that
\begin{equation}\label{eq:intro_Rosenthal_burkholder_ineq}
    \norm{\sup_{k\in\NB}\|\bM_k\|}_p\leq C_{1}\ \sqrt{p\vee \log d}\Big\|\norm{\qv{\bM}_\infty}^{1/2} \Big\|_p + C_{2}\ (p\vee \log d)\norm{ \sup_{k\in\NB}\|\bX_k\|}_p,
\end{equation}
for any $p\geq 2$.
Here, $\|\cdot\|$ is the spectral norm of the matrix, $\|\cdot\|_p$ is the $L^p$ norm of the random variable, and $C_{1}, C_{2}>0$ are universal constants.
While similar Rosenthal-type bounds for the spectral norm of matrix martingales were established in \citet{dirksen2023tuning}, our result provides a sharper dependence on the moment order $p$.
Specifically, this result matches the classic Hilbert space-valued Rosenthal-Burkholder inequality \citep[Theorem~4.1]{pinelis1994optimum}, with the only difference that an additional $\log d$ term appears. 
And the proof is based on the extrapolation method (Theorem~\ref{thm:good_lambda_inequality}) \citep{hitczenko1990best,hitczenko1994domination,pinelis1994optimum} and the techniques in the classic matrix concentration inequalities \citep{tropp2011freedmans, tropp2012user, tropp2015introduction, kroshnin2024bernstein}.

As a direct application of our Rosenthal-Burkholder inequality, we further establish a Burkholder-Davis-Gundy inequality for the spectral norm of continuous matrix local martingales $(\bM_t)_{t\geq 0}$ in Theorem~\ref{thm:martingale_BDG} by discretization, namely,
\begin{equation}\label{eq:intro_BDG}
    \begin{aligned}
       \norm{\sup_{0\le t\le \tau}\|\bM_t\|}_p\leq C_{1}\sqrt{p\vee \log d}\Big\|\norm{\qv{\bM}_\tau}^{1/2} \Big\|_p ,
    \end{aligned}
\end{equation}
for any $p\geq 2$ and stopping time $\tau$ such that $\|\norm{\qv{\bM}_\tau}^{1/2} \|_p<\infty$.

Regarding Markov chains, based on the Rosenthal-Burkholder inequality in Eqn.~\eqref{eq:intro_Rosenthal_burkholder_ineq}, we further derive Rosenthal inequalities for the spectral norm of sums of time-homogeneous matrix functions of ergodic Markov chains. Our proof combines the Poisson equation approach (see Eqn.~\eqref{eq:Poisson_equation}) \citep{durmus2023rosenthal} with a blocking argument and a matrix Bernstein inequality on a mixing-time skeleton chain \citep{neeman2024concentration}.
In particular, we consider a stationary Markov chain $(Z_k)_{k\in\NB}$ in a measurable space $\gZ$ with a Markov kernel $\rQ$ and a stationary distribution $\mu$, and a square-integrable self-adjoint matrix function
$\bF\in L^2_{\mu,0}(\HB^d)$ with $\mu(\bF^2) \coloneq \int_{\gZ}\bF^2(z)\mu(dz)< \infty$ and $\mu(\bF)=\bm{0}_{d\times d}$.
Our goal is to give an upper bound for the spectral norm of the summation $\bS_n \coloneq \sum_{k=0}^{n-1}{\bF}(Z_k)$.
By the Markov chain central limit theorem (CLT) \citep{jones2004markov}, under some regularity conditions, it holds that $n^{-1/2}\bS_n$ converges weakly to a random matrix with Gaussian entries and satisfies $\lim_{n\to\infty}n^{-1}\EB[\bS_n^2]=\bSigma_\mu(\bF)$, where $\bSigma_\mu(\bF)$ is the long-run variance defined as
\begin{equation}\label{eq:long_term_variance}
    \bSigma_\mu(\bF)\coloneq \mu\prn{{\bF}^2}+\sum_{k=1}^\infty\mu\prn{{\bF}\rQ^k{\bF}+ \prn{\rQ^k{\bF}}{\bF}},
\end{equation}
where $(\rQ \bF)(\cdot)\coloneq \int_{\gZ} \bF(z)\rQ(\cdot, dz)$.
Under the uniform geometric ergodicity assumption (Assumption~\ref{assumption:uniform_geo_ego}), Theorem~\ref{thm:markov_rosenthal} shows that, for any $p\geq2$ and $n\geq1$,
\begin{equation}\label{eq:intro_rosenthal_markov}
    \frac{1}{n}\norm{\sup_{k\leq n}\|\bS_k\|}_{p}
    \lesssim\sqrt{\frac{(p\vee\log d)\norm{\bSigma_\mu(\bF)}}{n}}
    +\frac{(p\vee\log d)\tmix}{n}\norm{\norm{\bF}}_\infty.
\end{equation}
The leading term depends on the long-run variance $\bSigma_\mu(\bF)$, matching the Markov chain CLT, while the residual term has the standard dependence on $\tmix$.
In the typical regime $\norm{\bSigma_\mu(\bF)}\simeq\tmix\|\|\bF\|\|_\infty^2$, the leading term dominates once $n\gtrsim(p\vee\log d)\tmix$. Theorem~\ref{thm:markov_rosenthal_geo_V} gives the corresponding long-run variance bound for geometrically $V$-ergodic chains and unbounded matrix functions.

By Markov's inequality, Corollary~\ref{cor:markov_hoeffding_uniform} and Corollary~\ref{cor:markov_bernstein_uniform} translate Eqn.~\eqref{eq:intro_rosenthal_markov} into Hoeffding-type and Bernstein-type high-probability upper bounds. For any $\delta\in(0,1)$, with probability at least $1-\delta$,
\begin{equation}\label{eq:intro_hoeffding_markov}
    \frac{1}{n}\sup_{k\leq n}\|\bS_k\|
    \lesssim\sqrt{\frac{\log(ed/\delta)\tmix\norm{\norm{\bF}}_{\infty}^2}{n}}
    +\frac{\log(ed/\delta)\tmix\norm{\norm{\bF}}_\infty}{n},
\end{equation}
\begin{equation}\label{eq:intro_bernstein_markov}
    \frac{1}{n}\sup_{k\leq n}\|\bS_k\|
    \lesssim\sqrt{\frac{\log(ed/\delta)\norm{\bSigma_\mu(\bF)}}{n}}
    +\frac{\log(ed/\delta)\tmix\norm{\norm{\bF}}_\infty}{n},
\end{equation}
where $\tmix\|\|\bF\|\|_\infty^2$ is a variance proxy that upper bounds $\|\bSigma_\mu(\bF)\|$ up to a universal constant.
\citet[Theorem~2.6]{neeman2024concentration} gave a Bernstein inequality with a leading term relying on the variance proxy $\gamma^{-1}\|\mu(\bF^2)\|$ which is also an upper bound of $\|\bSigma_\mu(\bF)\|$, where $\gamma$ is the absolute $L^2$-spectral gap.
Our leading term instead retains the long-run variance itself.
We apply our results in covariance estimation, PCA
with Markovian data, and obtain improved bounds in these applications.

The martingale inequalities \eqref{eq:intro_Rosenthal_burkholder_ineq} and \eqref{eq:intro_BDG} also admit dimension-free versions, where the ambient dimension $d$ is replaced by the effective rank (see Eqn.~\eqref{eq:effective_rank}) of an upper bound for the quadratic variation. These results extend to linear operators in infinite-dimensional Hilbert spaces.

Finally, we discuss other extensions of these results, which will not be elaborated on in this paper due to space constraints.
First, with some effort, we can prove using the same techniques that the inequality~\eqref{eq:intro_rosenthal_markov} also holds for geometrically ergodic Markov chains with respect to the Wasserstein semimetric (instead of total variation) \citep{hairer2011asymptotic} as in \citep{durmus2023rosenthal}.
This ergodicity assumption has applications in Markov chain Monte Carlo algorithms \citep{eberle2014error,hairer2014spectral}.
Second, based on the Burkholder-Davis-Gundy inequality~\eqref{eq:intro_BDG}, we can extend the results for Markov chains to continuous Markov processes.

\paragraph{Related Work}

Until now, probabilistic inequalities for the spectral norm of matrix martingales have been extensively studied. 
For example, \citet{oliveira2009concentration,tropp2011freedmans} proposed a Freedman inequality for the spectral norm of matrix martingales; \citet{minsker2017some,kroshnin2024bernstein} further derived dimension-free versions in different forms; and \citet{bacry2018concentration} established a Freedman inequality for the spectral norm of continuous-time matrix martingales.
Regarding Rosenthal-type moment inequalities under the spectral norm, \citet{chen2012masked} developed a sharp Rosenthal inequality in the independent case; \citet{jirak2024concentration} further extend the inequality to its dimension-free version.

Notably, \citet{dirksen2023tuning} also established two-sided Rosenthal-Burkholder inequalities for the spectral norm of matrix martingales. 
As noted in their work, their approach, which is rooted in the study of $L^p$-valued martingales \citep{dirksen2019q}, more generally yields moment inequalities for all Schatten $p$-norms for any $p > 1$. 
However, these results do not achieve the optimal constants for the spectral norm that we establish in Theorem~\ref{thm:martingale_rosenthal}. 
A detailed comparison of the constants and techniques is provided in the Discussion in Section~\ref{subsection:matrix_rosenthal_burkholder}.

Concurrently and independently,
\citet{maitre2025matrix} also established a Burkholder-Davis-Gundy inequality~\eqref{eq:intro_BDG} for the spectral norm of continuous matrix martingales using techniques from stochastic analysis.

Additionally, in the field of noncommutative probability \citep{pisier1997non,junge2003noncommutative,randrianantoanina2007conditioned,junge2008noncommutative,junge2013noncommutative,junge2015noncommutative,jiao2023noncommutative,malek2024noncommutative,jiao2025noncommutative}, the probabilistic inequalities for noncommutative martingales proposed therein can also be translated into probabilistic inequalities for the Schatten $p$-norms of matrix martingales. 
However, the existing results do not imply ours; see the Discussion in Section~\ref{subsection:matrix_rosenthal_burkholder} for details.

Regarding probability inequalities for Markov chains under the spectral norm, some existing work \citep{garg2018matrix,qiu2020matrix,neeman2024concentration,wu2025uncertainty,van2025matrix} has derived matrix concentration inequalities for Markov chains under various ergodicity assumptions.
Here, we elaborate on the differences between these existing results and ours in detail.

\citet{garg2018matrix,qiu2020matrix,neeman2024concentration,wu2025uncertainty} considered Markov chains with a non-zero absolute $L^2$-spectral gap (the gap between $1$ and the operator norm of $\rQ$ over the subspace $L^2_{\mu,0}$), using the multi-matrix Golden--Thompson inequality \citep{sutter2017multivariate,garg2018matrix} and the Le\'on--Perron argument \citep{leon2004optimal}.
These bounds are stated in terms of the spectral gap of the original kernel $\rQ$. Such an assumption implies geometric ergodicity with some Lyapunov function (Assumption~\ref{assumption:V_ego}) when the Markov chain is $\psi$-irreducible and aperiodic \citep[Theorem~1.3]{kontoyiannis2012geometric}, but the converse need not hold.
This assumption does not hold in many examples, especially for non-reversible Markov chains. 
For example, \citet[Example~A.4]{huang2024bernstein} provided a finite-state, non-reversible, irreducible, and aperiodic Markov chain; 
due to the finiteness of its state space, this Markov chain is uniformly geometrically ergodic (Assumption~\ref{assumption:uniform_geo_ego}), but its absolute $L^2$-spectral gap is zero. 
This limits their direct application to general non-reversible chains, which arise in reinforcement learning problems with Markovian data \citep{duan2022policy}. In contrast, uniform geometric ergodicity ensures that the skeleton kernel $\rQ^{2\tmix}$ has a constant absolute $L^2$-spectral gap (Lemma~\ref{lem:skeleton_gap}). We apply the Bernstein inequality of \citet{neeman2024concentration} only to this skeleton chain and combine it with the Poisson martingale decomposition. This preserves $\|\bSigma_\mu(\bF)\|$ in the leading term instead of replacing it by $\gamma^{-1}\|\mu(\bF^2)\|$.

In one recent work, \citet{van2025matrix} considered $\psi$-mixing Markov chains and employed the universality-based approach \citep{bandeira2023matrix,brailovskaya2024universality} in their proof.
Compared with the classic matrix concentration inequalities derived in our paper, the universality-based concentration inequalities they obtained are sharper in cases involving matrices with specific structures, especially when matrices in the summation are far from being commutative. 
In their upper bounds, the coefficient of the long-run variance $\|\bSigma_\mu(\bF)\|$ no longer includes $\log d$. 
Instead, $\log d$ appears in a new variance term, which is a lower bound of $\|\bSigma_\mu(\bF)\|$. 
For a detailed introduction of universality-based concentration inequalities, one may refer to the seminal work \citep{bandeira2023matrix,brailovskaya2024universality}.
However, the results in \citep{van2025matrix} have the following limitations. 
On the one hand, $\psi$-mixing is a stronger ergodicity assumption than a non-zero absolute $L^2$-spectral gap and uniform geometric ergodicity (Assumption~\ref{assumption:uniform_geo_ego}). 
Markov chains with infinite states generally fail to satisfy this condition, which limits the scope of application of the results.
On the other hand, the approach they employed remains difficult to extend to infinite-dimensional Hilbert spaces.

\paragraph{Notation}
Here and later,  
$a\wedge b\coloneq \min\{a,b\}$, $a\vee b\coloneq  \max\{a,b\}$ for any $a,b\in\RB$.
We also denote $\RB_+\coloneq [0,\infty)$, $\NB^*\coloneq \{1,2,3,\ldots\}$.
The symbol ``$\lesssim$'' (resp. ``$\gtrsim$'') means no larger (resp. smaller) than up to a multiplicative universal constant, and $a\simeq b$ means $a\lesssim b$ and $a\gtrsim b$ hold simultaneously.
We denote by $\ind$ the indicator function and $\abs{\xi}$ the total variation measure of a signed measure $\xi$.
To simplify the presentation, we adopt the shorthand $\mathbb{E}^{\alpha}[\cdot]$ to represent $(\mathbb{E}[\cdot])^{\alpha}$ for any $\alpha>0$.
For any $p\geq 1$ and a real-valued random variable $X$, we denote by $\|X\|_p=\EB^{1/p}[|X|^p]$ the $L^p$ norm of $X$.

We denote by $\bI_d\in\RB^{d{\times} d}$ the $d\times d$ identity matrix, $\bm{0}_{d_1{\times} d_2}\in\RB^{d_1{\times} d_2}$ the $d_1{\times} d_2$ zero matrix,
$\norm{\bu}$ the Euclidean norm of vector $\bu$,
$\norm{\bB}=\sup_{\bu\colon\|\bu\|\leq 1}\|\bB \bu\|$ the spectral norm of matrix $\bB$, $\bB^*$ the conjugate transpose of matrix $\bB$, and $[\bB]_{ij}$ the $(i,j)$-entry of matrix $\bB$.
For any self-adjoint matrix $\bB\in\HB^d$ ($\bB=\bB^*$) with spectral decomposition $\bB=\sum_{i=1}^d \lambda_i\bu_i\bu_i^*$, we denote by $\tr{\bB}=\sum_{i=1}^d \lambda_i$ the trace of matrix $\bB$, and we define $|\bB|\coloneq \sum_{i=1}^d |\lambda_i|\bu_i\bu_i^*$ , $\bB\wedge a\coloneq \sum_{i=1}^d (\lambda_i \wedge a)\bu_i\bu_i^*$ for any $a\in\RB$.
For any $\bm B_1, \bm B_2\in\HB^d$, $\bm B_1\preccurlyeq\bm B_2$ stands for $\bm B_2{-}\bm B_1$ is positive semi-definite (PSD),
and we denote by $\HB^d_+$ the set of all $d\times d$ PSD matrices.

\paragraph{Organization}
The remainder of this paper is organized as follows. 
In Section~\ref{Section:main_results}, we present our main results.
In Section~\ref{Section:application}, we apply our results in covariance estimation, PCA problems with Markovian data, and obtain improved error bounds.
The proof is outlined in Section~\ref{Section:proof_outline}.
In Section~\ref{Section:conclusion}, we conclude our work.
Details of the proof are given in the appendices.

%% file: main_results_martingale.tex
In this section, we present our main results: the Rosenthal-Burkholder inequality for the spectral norm of discrete-time matrix martingales and Burkholder-Davis-Gundy inequality for the spectral norm of continuous matrix local martingales in Section~\ref{subsection:matrix_rosenthal_burkholder}, as well as matrix Rosenthal, Hoeffding, and Bernstein inequalities under the spectral norm for ergodic Markov chains in Section~\ref{subsection:markov_matrix_rosenthal}. 
We will also give  dimension-free versions of the martingale results.

We find that the dimension-free results can be applied to linear operator-valued martingales in infinite-dimensional separable Hilbert spaces. 
For brevity, we will not change the notation, and our dimension-free results still hold if $\HB^d$ is replaced by the space of all self-adjoint operators in the separable Hilbert space $\gX$, and $\CB^{d_1{\times}d_2}$ is replaced by the space of all bounded linear operators from the separable Hilbert space $\gX_1$ to the separable Hilbert space $\gX_2$.

\subsection{The Rosenthal-Burkholder Inequality under the Spectral Norm}\label{subsection:matrix_rosenthal_burkholder}

Let $(\Omega,\sF,\PB)$ be a probability space and $\{\emptyset,\Omega\}=\sF_{-1}\subseteq\sF_0\subseteq \sF_1\subseteq\ldots\subseteq \sF$ be the filtration.
We denote by $\EB_k$ the expectation conditioned on $\sF_k$.
We consider the $d$-dimensional self-adjoint matrix martingale $(\bM_k)_{k\in\NB}$ that is adapted to the filtration $(\sF_k)_{k\in\NB}$.
The matrix martingale $(\bM_k)_{k\in\NB}$ satisfies 
\begin{equation*}
    \EB_{k-1}\bM_k=\bM_{k-1}\quad\text{  and  }\quad \EB\norm{\bM_k}<\infty\quad\text{ for }\ k\in\NB.
\end{equation*}
We denote by $(\bX_k)_{k\in\NB}$ the martingale difference sequence, that is, $\bX_k=\bM_k-\bM_{k-1}\in\HB^d$ ($\bM_{-1}\coloneq \bm{0}_{d{\times}d}$).
And we denote by $(\qv{\bM}_k)_{k\in\NB}$ the predictable quadratic variation sequence, that is, $\qv{\bM}_k=\sum_{j=0}^k \EB_{j-1}\bX_j^2\in\HB^d_+$.
Throughout this section, we assume that $\qv{\bM}_\infty\coloneq \lim_{k\to\infty}\qv{\bM}_k<\infty$ element-wisely almost surely.
As a result, by the martingale convergence theorem \citep[Theorem~4.5.2]{durrett2019probability}, $\lim_{k\to\infty}{\bM}_k=\bM_\infty\in\HB^d$ almost surely.

\textit{Rectangular Case: }For a $d_1{\times} d_2$ rectangular matrix martingale $(\bW_k)_{k\in\NB}$ adapted to $(\sF_k)_{k\in\NB}$, we denote by $\bD_k=\bW_k-\bW_{k-1}$ the martingale difference and $\|\langle \bW\rangle_\infty\|\coloneq \|\sum_{k=0}^\infty\EB_{k-1}[\bD_{k}\bD_{k}^{*}]\|\vee\|\sum_{k=0}^\infty\EB_{k-1}[\bD_{k}^{*}\bD_{k}] \|$ the proxy of the quadratic variation.
In this case, we also assume that
\[
\sum_{k=0}^\infty\EB_{k-1}[\bD_{k}\bD_{k}^{*}]<\infty \mbox{ and } \sum_{k=0}^\infty\EB_{k-1}[\bD_{k}^{*}\bD_{k}]<\infty \quad \mbox{element-wisely almost surely}.
\]

\textit{Dimension-Free: }As mentioned in \citet{minsker2017some}, in the case of ultra-high dimensions or even infinite dimensions, it is desirable that the results do not depend on the dimension $d$ of the ambient space. 
To this end, we resort to the notion of effective rank \citep{koltchinskii2016asymptotics}.
For any $\bUp\in\HB^d_+$, its effective rank is defined as
\begin{equation}\label{eq:effective_rank}
    r(\bUp)=\frac{\tr(\bUp)}{\norm{\bUp}}.
\end{equation}
It is easily seen that $r(\bUp)\leq d$ and $r(\bUp)$ can be much smaller than $d$ if the eigenvalues of $\bUp$ rapidly decay to $0$.

Now, we are ready to state one of the main results, the Rosenthal-Burkholder inequality for the spectral norm of matrix martingales in the following theorem.
\begin{theorem}[Rosenthal-Burkholder Inequality]\label{thm:martingale_rosenthal}
For any $p\geq 2$, suppose that
\begin{equation*}
    \quad \Big\|\norm{\qv{\bM}_\infty}^{1/2} \Big\|_p<\infty \quad \text{ and }\quad\norm{ \sup_{k\in\NB}\|\bX_k\|}_p<\infty.
\end{equation*}
Then it holds that
\begin{equation}\label{eq:used_burkholder_ineq}
    \begin{aligned}
       \norm{\sup_{k\in\NB}\|\bM_k\|}_p\leq C_{1}\sqrt{p\vee \log d}\Big\|\norm{\qv{\bM}_\infty}^{1/2} \Big\|_p + C_{2}(p\vee\log d)\norm{ \sup_{k\in\NB}\|\bX_k\|}_p,
    \end{aligned}
\end{equation}
where $C_{1}=87$, $C_{2}=50$. 
When $p\geq 117$, the constants can be smaller: $C_{1}=64$, $C_{2}=28$.  

Furthermore, the following extensions hold:
\begin{enumerate}[(a)]
    \item \textit{Rectangular Matrix:}
The above inequality also holds for the rectangular matrix martingale $(\bW_k)_{k\in\NB}$, in which we replace $\bM$ by $\bW$, $\bX$ by $\bD$, and $d$ by $d_1+d_2$.
\item \textit{Dimension-Free:} In the self-adjoint case, if we further assume that there exists a $\bUp\in\HB^d_+$ such that $\langle\bM\rangle_\infty\preccurlyeq \bUp$ almost surely, then the above inequality also holds when we replace $d$ by $e\cdot r(\bUp\wedge \|\|\qv{\bM}_\infty\|\|_{\frac{p}{2}})$, and double the right-hand side.

In the rectangular case (a), if we further assume that there exist $\bUp_1\in\HB_+^{d_1}$, $\bUp_2\in\HB_+^{d_2}$ such that $\sum_{k=0}^\infty\EB_{k-1}[\bD_{k}\bD_{k}^{*}]\preccurlyeq\bUp_1$ and  $\sum_{k=0}^\infty\EB_{k-1}[\bD_{k}^{*}\bD_{k}]\preccurlyeq\bUp_2$ almost surely, then the above inequality also holds, while we replace $\bM$ by $\bW$, $\bX$ by $\bD$, $d$ by $2[r(\bUp_1\wedge \|\|\qv{\bW}_\infty\|\|_{\frac{p}{2}})\vee r(\bUp_2\wedge \|\|\qv{\bW}_\infty\|\|_{\frac{p}{2}})]$, and double the right-hand side. 
\end{enumerate}
\end{theorem}

By localization, our Rosenthal-Burkholder inequalities also hold if we consider matrix local martingales and replace $k\in\NB$ by $0\leq k\leq \tau$ for any stopping time $\tau$; we omit the details for brevity.
We also present another version of the Rosenthal-Burkholder inequality: $\|\sup_{k\in\NB}\|\bM_k\|\|_p\lesssim [(p/\log p)\vee \log d]( \|\norm{\qv{\bM}_\infty}^{1/2} \|_p  + \| \sup_{k\in\NB}\|\bX_k\|\|_p)$, which can be found in Section~\ref{Subsection:proof_outline_burkholder}.

For self-adjoint matrix martingales $(\bM_k)_{k\in\NB}$, the proof strategy is to combine the extrapolation method (\textit{good $\lambda$ inequality}) \citep{hitczenko1990best,hitczenko1994domination,pinelis1994optimum}
with techniques in matrix concentration inequalities (especially matrix Bennett inequalities) \citep{tropp2011freedmans, tropp2012user, tropp2015introduction, kroshnin2024bernstein}.
A detailed outline
of the proof can be found in Section~\ref{Subsection:proof_outline_burkholder}.

For the rectangular matrix martingale $(\bW_k)_{k\in\NB}$, we can directly apply the previous results to the self-adjoint matrix martingale $(\sH(\bW_k))_{k\in\NB}$ to obtain the desired inequality, where $\sH$ is the Hermitian dilation operator \citep{paulsen2002completely}.
For a rectangular matrix $\bB$, its Hermitian dilation $\sH\prn{\bB}$ is a self-adjoint matrix, which satisfies
\begin{equation*}
    \sH\prn{\bB}\coloneq \begin{bmatrix} 
\bm{0} & \bB \\
\bB^* & \bm{0}
\end{bmatrix}, \quad \sH\prn{\bB}^2=\begin{bmatrix} 
\bB\bB^* & \bm{0} \\
\bm{0} & \bB^*\bB
\end{bmatrix},\quad \norm{\sH\prn{\bB}}=\norm{\bB}.
\end{equation*}

We note that doubling the r.h.s. in Theorem~\ref{thm:martingale_rosenthal} (b) is not essential; we could reduce this multiplier from $2$ to a constant close to $1$ by increasing the effective rank term that only appears as the logarithmic term. 
However, for brevity, we use this presentation and do not strive for the optimal constants.

A direct application of the Rosenthal-Burkholder inequality (Theorem~\ref{thm:martingale_rosenthal}) is to show the Burkholder-Davis-Gundy inequality for the spectral norm of continuous matrix local martingales by discretization in the following theorem.
\begin{theorem}[Burkholder-Davis-Gundy Inequality]\label{thm:martingale_BDG}
Let $(\bM_t)_{t\geq 0}$ be a continuous local martingale taking values in $\HB^d$. We denote by $(\qv{\bM}_t)_{t\geq 0}$ its quadratic variation process, that is, $[\qv{\bM}_t]_{ij}=\sum_{k=1}^d\qv{[\bM]_{ik},[\bM]_{kj}}_t$. For any $p\geq 2$ and stopping time $\tau$, suppose that  $\|\norm{\qv{\bM}_\tau}^{1/2} \|_p<\infty$. Then it holds that
\begin{equation*}
    \begin{aligned}
       \norm{\sup_{0\le t\le \tau}\|\bM_t\|}_p\leq C_{1}\sqrt{p\vee \log d}\Big\|\norm{\qv{\bM}_\tau}^{1/2} \Big\|_p ,
    \end{aligned}
\end{equation*}
where $C_{1}$ is defined in Theorem~\ref{thm:martingale_rosenthal}.  
\end{theorem}
The proof can be found in Appendix~\ref{appendix:proof_bdg_inequality}.

\paragraph{Discussion}
In the independent case, \citet[Theorem~A.1]{chen2012masked} provided a Rosenthal inequality for the spectral norm,
and \citet[Theorem~3]{jirak2024concentration} further gave a dimension-free version.
We extend their results from the independent case to the martingale case.

In the martingale setting, \citet[Theorem~6]{dirksen2023tuning} established two-sided Rosenthal-Burkholder inequalities for the spectral norm of matrix martingales. 
However, the constants in their inequalities involve a suboptimal dependence on the moment order $p$. 
In contrast, our Theorem~\ref{thm:martingale_rosenthal} achieves a sharper dependence on $p$ with the scaling of $\sqrt{p}$ and $p$, matching the optimal constants in the classic scalar case \citep{pinelis1985estimates}. 
While a common element in our proof and that of \citet[Theorem~6]{dirksen2023tuning} is the use of decoupling of martingale difference sequences, we combine the extrapolation method (\textit{good $\lambda$ inequality}) \citep{hitczenko1990best,hitczenko1994domination,pinelis1994optimum}
with techniques in matrix concentration inequalities \citep{tropp2011freedmans, tropp2012user, tropp2015introduction, kroshnin2024bernstein}. 
This leads to an improvement of the constants that features the same order as those established in the case of independent martingale differences.

In the field of noncommutative probability, \citet{pisier1997non,junge2003noncommutative,randrianantoanina2007conditioned,junge2008noncommutative,junge2013noncommutative,junge2015noncommutative,jiao2023noncommutative,malek2024noncommutative,jiao2025noncommutative} investigated probability inequalities for noncommutative martingales.
In our random matrix setting, their results can be translated into probability inequalities concerning the Schatten $p$-norm, rather than the spectral norm that is the focus of this paper.
In particular, in the independent case, \citet[Theorem~0.4]{junge2013noncommutative} gives the following Rosenthal inequality
\begin{equation*}
    \begin{aligned}
       \EB^{{1}/{p}}\brk{\tr\prn{ \abs{\bM_{\infty}}^p}}\lesssim \sqrt{p}\EB^{{1}/{p}}\brk{\tr\prn{\qv{\bM}_\infty^{p/2}}}  + p\EB^{{1}/{p}}\brk{\sup_{k\in\NB}\tr\prn{ \abs{\bX_{k}}^p}},
    \end{aligned}
\end{equation*}
for $p\geq 2.5$.
However, in the martingale case, whether a Rosenthal-Burkholder inequality of the following form holds remains an open problem \citep[Problem~2.5]{junge2015noncommutative}
\begin{equation*}
    \begin{aligned}
       \EB^{{1}/{p}}\brk{\tr\prn{ \abs{\bM_{\infty}}^p}}\lesssim \sqrt{p}\EB^{{1}/{p}}\brk{\tr\prn{\qv{\bM}_\infty^{p/2}}}  + f(p)\EB^{{1}/{p}}\brk{\sum_{k\in\NB}\tr\prn{ \abs{\bX_{k}}^p}},
    \end{aligned}
\end{equation*}
where the coefficient of the dominant term (quadratic variation) has an optimal $\sqrt{p}$ dependence on $p$, and $f(p)$ is some function for which they conjectured $f(p)=p$ as in the commutative case.

%% file: main_results_markov.tex
\subsection{Probability Inequalities under the Spectral Norm for Markov Chains}\label{subsection:markov_matrix_rosenthal}

Consider the Markov kernel $\rQ\colon \gZ\times \sZ\to [0,1]$ on the measurable space $(\gZ,\sZ)$ with a unique stationary distribution $\mu\in\gP$, where $\gP$ is the set of all probability distributions on $(\gZ,\sZ)$.
For any $\xi\in\gP$, we denote by $\PB_\xi$ (resp. $\EB_\xi$) the probability (resp. expectation) of the canonical process $(Z_k)_{k\in\NB}$ with initial distribution $\xi$, and by $(\sF_k)_{k\in\NB}$ the natural filtration, where $\sF_k=\sigma(Z_0,\ldots, Z_k)$.
That is, $(Z_k)_{k\in\NB}$ is a stochastic process in $(\gZ^\NB,\sZ^\NB,\PB_\xi)$ adapted to $(\sF_k)_{k\in\NB}$, which satisfies $\PB_\xi(Z_0\in A)=\xi(A)$, $\PB_\xi(Z_{k+1}\in A\mid \sF_{k})=\PB_\xi(Z_{k+1}\in A\mid Z_{k})=\rQ(Z_{k}, A)$ $\PB_\xi$-almost surely for any $A\in\sZ$ and $k\in\NB$.
Thus, $(Z_k)_{k\in\NB}$ is a Markov chain.

For any $p\geq1$, $\xi\in\gP$ and random variable $X$ defined on $(\gZ^\NB,\sZ^\NB,\PB_\xi)$, we denote by $\|X\|_{\xi,p}=(\EB_{\xi}[|X|^p])^{1/p}$ the $L^p$ norm of $X$ under the initial distribution $\xi$.

For any measurable function $f$ in $(\gZ,\sZ)$ and distribution $\xi\in\gP$, we denote $\norm{f}_{\infty}=\sup_{z\in\gZ}\abs{f(z)}$ ($f$ is real-valued), $\xi(f)=\int_\gZ f(z) \xi(dz)$, $(\xi \rQ)(\cdot)=\int_{\gZ}\rQ(z,\cdot)\xi(dz)$, and $(\rQ f)(\cdot)=\int_{\gZ} f(z)\rQ(\cdot, dz)$.

We now introduce the ergodicity assumptions considered in this paper. In this section, we focus on uniformly geometrically ergodic Markov chains; the more general geometrically $V$-ergodic chains are considered in Section~\ref{subsection:geo_V_ergo}.
We recall that an irreducible and aperiodic Markov chain with a finite state space is always uniformly geometrically ergodic.

\begin{assumption}[Uniformly Geometrically Ergodic Markov Chains]\label{assumption:uniform_geo_ego}
    There exists a $\tmix\in\NB$ such that for any $k\in\NB$,
\begin{equation*}
    \sup_{z\in\gZ}\norm{\rQ^k(z,\cdot)-\mu}_{\TV}\leq 2 \prn{1/4}^{\lfloor k/\tmix\rfloor},
\end{equation*}
where $\norm{\xi}_{\TV}\coloneq \abs{\xi}(\gZ)$ is the total variation norm of the signed measure $\xi$ on $(\gZ,\sZ)$.
\end{assumption}

We consider a measurable self-adjoint matrix function $\bF\colon \gZ\to\HB^d$.
We assume $\mu(\bF)=\int_{\gZ}\bF(z)\mu(dz)=\bm{0}_{d{\times}d}$ and $\|\|\bF \|\|_\infty=\sup_{z\in\gZ}\|\bF(z)\|<\infty$.
Our goal is to give counterparts to Rosenthal, Hoeffding, and Bernstein inequalities for the spectral norm of the summation $\bS_n=\sum_{k=0}^{n-1}{\bF}(Z_k)$.
We observe that the boundedness condition is indispensable. 
As argued by \citet[Theorem~4]{fan2021hoeffding}, even for scalar unbounded functions, it is impossible to obtain Hoeffding-type inequalities for $\bS_n$ under Assumption~\ref{assumption:uniform_geo_ego} or a non-zero absolute $L^2$-spectral gap, whereas our Rosenthal inequality implies the Hoeffding-type inequality for $\bS_n$.

Recall that the Markov chain CLT \citep[Corollary~5]{jones2004markov} says that $n^{-1/2}\bS_n$ converges weakly to a Gaussian distribution and satisfies $\lim_{n\to\infty}n^{-1}\EB_{\xi}[\bS_n^2]=\bSigma_\mu(\bF)$ for any $\xi\in\gP$, where $\bSigma_\mu(\bF)$ is the long-run variance defined as
\begin{equation*}
    \bSigma_\mu(\bF)\coloneq \mu\prn{{\bF}^2}+\sum_{k=1}^\infty\mu\prn{{\bF}\rQ^k{\bF}+ \prn{\rQ^k{\bF}}{\bF}}.
\end{equation*}
We aim to obtain a Rosenthal inequality whose leading term is given by the long-run variance $\bSigma_\mu(\bF)$, matching the Markov chain CLT.

\textit{Rectangular Case: }In this case, we consider a rectangular matrix function $\bR\colon\gZ\to\CB^{d_1{\times}d_2}$ that satisfies $\mu(\bR)=\int_{\gZ}\bR(z)\mu(dz)=\bm{0}_{d_1{\times}d_2}$. We denote by $\bV_n=\sum_{k=0}^{n-1}\bR(Z_k)$ the summation and
\begin{equation}\label{eq:long_term_variance_rect}
    \norm{\bSigma_\mu(\bR)}\coloneq \norm{\mu\prn{\bR\bR^*}{+}\sum_{k=1}^\infty\mu\prn{{\bR}\rQ^k{\bR^*}{+} \prn{\rQ^k{\bR}}{\bR^*}}}\vee\norm{\mu\prn{\bR^*\bR}{+}\sum_{k=1}^\infty\mu\prn{{\bR^*}\rQ^k{\bR}{+} \prn{\rQ^k{\bR^*}}{\bR}}},
\end{equation}
called the proxy of the long-run variance.

Now we state one of our main results, the Rosenthal-type inequality under the spectral norm for uniformly geometrically ergodic Markov chains.
\begin{theorem}[Rosenthal Inequality for Uniformly Geometrically Ergodic Markov Chains]\label{thm:markov_rosenthal}
Under Assumption~\ref{assumption:uniform_geo_ego} and $\|\|\bF\|\|_\infty<\infty$, there exist universal constants $D_1,D_2,D_3>0$ such that, for any $p\geq 2$, $n\geq 1$, and $\xi\in\gP$,
\begin{equation}\label{eq:sharp_stationary_rosenthal}
    \frac{1}{n}\norm{\sup_{k\leq n}\norm{\bS_k}}_{\mu,p}
    \leq D_1\sqrt{\frac{(p\vee\log d)\norm{\bSigma_\mu(\bF)}}{n}}
    +D_2\frac{(p\vee\log d)\tmix}{n}\norm{\norm{\bF}}_{\infty},
\end{equation}
and
\begin{equation}\label{eq:from_non_stationary_to_stationary}
    \frac{1}{n}\norm{\sup_{k\leq n}\norm{\bS_k}}_{\xi,p}
    \leq \frac{1}{n}\norm{\sup_{k\leq n}\norm{\bS_k}}_{\mu,p}
    +D_3\frac{p\tmix}{n}\norm{\norm{\bF}}_{\infty}.
\end{equation}
Consequently, for another universal constant $D_4>0$,
\begin{equation}\label{eq:sharp_general_rosenthal}
    \frac{1}{n}\norm{\sup_{k\leq n}\norm{\bS_k}}_{\xi,p}
    \leq D_1\sqrt{\frac{(p\vee\log d)\norm{\bSigma_\mu(\bF)}}{n}}
    +D_4\frac{(p\vee\log d)\tmix}{n}\norm{\norm{\bF}}_{\infty}.
\end{equation}
The same inequalities hold in the rectangular case if $\|\|\bR\|\|_\infty<\infty$ and we replace $\bS$ by $\bV$, $\bF$ by $\bR$, $d$ by $d_1+d_2$, and $\norm{\bSigma_\mu(\bF)}$ by the proxy $\norm{\bSigma_\mu(\bR)}$ defined in Eqn.~\eqref{eq:long_term_variance_rect}.
\end{theorem}

The proof combines the Poisson equation $\bG-(\rQ\bG)=\bF$ with a blocking argument. The Poisson equation identifies the long-run variance as the expected predictable quadratic variation of a matrix martingale. 
We then apply a matrix Bernstein inequality from \citet{neeman2024concentration} to the conditional variances along the $2\tmix$-skeleton chain and use Theorem~\ref{thm:martingale_rosenthal} at the block endpoints. This gives the long-run variance in the leading term and the standard dependence $\tmix/n$ in the residual term. The proof is given in Section~\ref{Subsection:proof_outline_markov}.

By Assumption~\ref{assumption:uniform_geo_ego}, $\norm{\bSigma_\mu(\bF)}\leq (19/3)\tmix\norm{\norm{\bF}}_{\infty}^2$. Thus, Theorem~\ref{thm:markov_rosenthal} and Markov's inequality immediately yield the following Hoeffding and Bernstein inequalities.
\begin{corollary}[Hoeffding Inequalities for Uniformly Geometrically Ergodic Markov Chains]\label{cor:markov_hoeffding_uniform}
Under Assumption~\ref{assumption:uniform_geo_ego} and $\|\|\bF \|\|_\infty<\infty$,
for any $\delta\in(0, 1)$, $n\geq 1$, and $\xi\in\gP$, with $\PB_\xi$-probability at least $1-\delta$, it holds that
\begin{equation*}
   \frac{1}{n}\sup_{k\leq n}\norm{\bS_k}
   \leq D_5\sqrt{\frac{\log(ed/\delta)\tmix\norm{\norm{\bF}}_{\infty}^2}{n}}
   +D_6\frac{\log(ed/\delta)\tmix\norm{\norm{\bF}}_{\infty}}{n},
\end{equation*}
where $D_5,D_6>0$ are universal constants.
\end{corollary}

\begin{corollary}[Bernstein Inequalities for Uniformly Geometrically Ergodic Markov Chains]\label{cor:markov_bernstein_uniform}
Under Assumption~\ref{assumption:uniform_geo_ego} and $\|\|\bF \|\|_\infty<\infty$,
for any $\delta\in(0, 1)$, $n\geq 1$, and $\xi\in\gP$, with $\PB_\xi$-probability at least $1-\delta$, it holds that
\begin{equation*}
   \frac{1}{n}\sup_{k\leq n}\norm{\bS_k}
   \leq D_5\sqrt{\frac{\log(ed/\delta)\norm{\bSigma_\mu(\bF)}}{n}}
   +D_6\frac{\log(ed/\delta)\tmix\norm{\norm{\bF}}_{\infty}}{n}.
\end{equation*}
\end{corollary}
The rectangular counterparts follow from the same Hermitian dilation argument as in Theorem~\ref{thm:markov_rosenthal}.

%% file: main_results_exteinsion_V_ergo.tex
\subsubsection{Results for Geometrically \texorpdfstring{$V$}{V}-Ergodic Markov Chains}\label{subsection:geo_V_ergo}
The assumptions of uniform geometric ergodicity (Assumption~\ref{assumption:uniform_geo_ego}) and bounded functions are excessively strong and fail to hold in many applications. 
We now introduce the more general geometrically ergodic Markov chains with Lyapunov function $V$ and present Rosenthal inequalities under the ergodicity assumption.
Let $V\colon \gZ\to[e,\infty)$ be a Lyapunov function satisfying $\mu(V)<\infty$. 
We denote by $\gP_V$ the set of all probability distributions whose elements $\xi $ satisfy $\xi(V)<\infty$.
For any real-valued measurable function $f$, we also define its $V$-norm as $\|f\|_V=\sup_{z\in\gZ}{\abs{f(z)}}/{V(z)}$.
We will see that, under the ergodicity assumption, for matrix functions $\bF$ whose growth rate does not exceed $V^{1/p}$, the convergence of $n^{-1}\bS_n$ in the $L^p$ norm can be obtained, which allows us to extend the result to unbounded functions.
In this section we omit the results in the rectangular case for brevity.
\begin{assumption}[Geometrically $V$-Ergodic Markov Chains]\label{assumption:V_ego}
    There exist $\tmix\in\NB$ and $\varkappa\geq 1$ such that for any $k\in\NB$,
\begin{equation*} 
    \sup_{z,z^\prime\in\gZ}\frac{\norm{\rQ^k(z,\cdot)-\rQ^k(z^\prime,\cdot)}_V}{V(z)+V(z^\prime)}\leq \varkappa\prn{\frac{1}{4}}^{\lfloor k/\tmix\rfloor},
\end{equation*}
where $\norm{\xi}_{V}\coloneq \abs{\xi}(V)$ is the
the $V$-total variation norm of the signed measure $\xi$ on $(\gZ,\sZ)$.
\end{assumption}
This assumption can be verified via the Foster-Lyapunov drift condition \citep{meyn2012markov}.
Assumption~\ref{assumption:V_ego} implies
that for any $\xi\in\gP_{V}$ and $k\in\NB$, 
\begin{equation*}
    \norm{\xi \rQ^k -\mu}_V\leq \varkappa\prn{\xi(V)+\mu(V)}\prn{1/4}^{\lfloor k/\tmix\rfloor}.
\end{equation*}
Hence, Assumption~\ref{assumption:uniform_geo_ego} is a special case of Assumption~\ref{assumption:V_ego} with $V\equiv e$ and $\varkappa=1$.

Now we are ready to present Rosenthal-type inequality under the spectral norm for geometrically $V$-ergodic Markov chains.
We first present the  result parallel to Theorem~\ref{thm:markov_rosenthal} whose leading term matches the Markov chain CLT, with the proof in Section~\ref{Subsection:proof_outline_markov}.
\begin{theorem}[Rosenthal Inequality for Geometrically $V$-Ergodic Markov Chains]\label{thm:markov_rosenthal_geo_V}
Under Assumption~\ref{assumption:V_ego} and $\|\|\bF \|\|_{V^{1/p}}<\infty$ for $p\geq 2$, we have that for any $n\geq 2$, and $\xi\in\gP_V$, 
\begin{equation*}
    \begin{aligned}
        \frac{1}{n}\norm{\sup_{k\leq n}\norm{\bS_k}}_{\mu,p}\leq& C_{1}\sqrt{\frac{(p\vee\log d)\norm{\bSigma_\mu(\bF)}}{n}}+D_{4}p^{3/2}\prn{p\vee\log d}^{3/4}\tmix^{5/4}n^{-3/4}\brk{\varkappa\mu(V)}^{2/p}\norm{\norm{\bF}}_{V^{1/p}}\\
        &+D_5p^{3/2}\prn{p\vee\log d}\tmix^{3/2}n^{-(1-1/p)}\brk{\varkappa\mu(V)}^{2/p}\norm{\norm{\bF}}_{V^{1/p}},
    \end{aligned}
\end{equation*}
\begin{equation*}
    \begin{aligned}
        \frac{1}{n}\norm{\sup_{k\leq n}\norm{\bS_k}}_{\xi,p}\leq \frac{1}{n}\norm{\sup_{k\leq n}\norm{\bS_k}}_{\mu,p}+\frac{32}{15n}\varkappa^{1/p}\prn{\xi(V)+\mu(V)}^{1/p}p\tmix\norm{\norm{\bF}}_{V^{1/p}},
    \end{aligned}
    \end{equation*}
    where $C_{1}$ is defined in Theorem~\ref{thm:martingale_rosenthal} and $D_4, D_5>0$ are universal constants, which can be traced in the proof.
\end{theorem}
The proof uses the same Poisson martingale decomposition as in the uniformly geometrically ergodic case, together with $V$-norm estimates that allow $\bF$ to be unbounded.

For comparison, the blocking argument in Lemma~\ref{lem:crude_rosenthal_markov} (b) gives the following bound in terms of the variance proxy $\tmix\|\|\bF\|\|_{V^{1/p}}^2$.
\begin{theorem}\label{thm:markov_crude_rosenthal_V}
Under Assumption~\ref{assumption:V_ego} and $\|\|\bF\|\|_{V^{1/p}}<\infty$ for $p\geq 2$, we have that, 
    \begin{equation*}
    \begin{aligned}
        \frac{1}{n}\norm{\sup_{k\leq n}\norm{\bS_k}}_{\mu,p}\leq& \frac{8}{3}C_{1}p\sqrt{\frac{(p\vee \log d)\tmix}{n}}\brk{\varkappa\mu(V)}^{1/p}\norm{\norm{\bF}}_{V^{1/p}} \\
        &+ \prn{\frac{16}{3}C_{2}(p\vee\log d)+\frac{11}{3}}\frac{p\tmix}{ n^{1-1/p}}\brk{\varkappa\mu(V)}^{1/p} \norm{\norm{\bF}}_{V^{1/p}}.
    \end{aligned}
    \end{equation*}
\end{theorem}
We can also use Markov inequality to obtain high-probability upper bounds similar to Corollary~\ref{cor:markov_hoeffding_uniform} and Corollary~\ref{cor:markov_bernstein_uniform}. 
For brevity, we omit them here.

\paragraph{Discussion} 
As discussed in Related Work, \citet{garg2018matrix,qiu2020matrix,neeman2024concentration,wu2025uncertainty}  derived matrix concentration inequalities for Markov chains with a non-zero absolute $L^2$-spectral gap; and \citet{van2025matrix} derived universality-based concentration inequalities for $\psi$-mixing Markov chains.
Our proof starts from the Poisson equation approach \citep{durmus2023rosenthal}. For the sharp uniformly geometrically ergodic result, we additionally apply the matrix Bernstein inequality of \citet{neeman2024concentration} to a skeleton chain whose absolute spectral gap is bounded below by a universal constant.
First, our results can be applied to general non-reversible geometrically ergodic Markov chains without assuming a non-zero absolute spectral gap or $\psi$-mixing condition.
Second, we do not need to assume the boundedness of functions under the geometric ergodicity assumption. 
Third, the leading term of our results matches the Markov chain CLT, rather than relying on a sub-optimal variance proxy as in \citep{garg2018matrix,qiu2020matrix,neeman2024concentration,wu2025uncertainty}.

For uniformly geometrically ergodic chains and bounded functions, Theorem~\ref{thm:markov_rosenthal} has the sharp $\tmix/n$ residual term. Extending this dependence to the geometrically $V$-ergodic and unbounded setting of Theorem~\ref{thm:markov_rosenthal_geo_V} remains open.
The Poisson equation approach is also restricted to time-homogeneous functions (additive functionals of Markov chains), whereas the results in \citet{neeman2024concentration} apply to time-inhomogeneous functions. 
This makes our results inapplicable to problems such as stochastic approximation with Markovian data \citep{li2023online,mou2024optimal,huo2024effectiveness,samsonov2025statistical}.
A final direction is suggested by the universality-based inequalities in \citep{van2025matrix}, which yield sharper bounds for some structured matrices under the stronger $\psi$-mixing condition and a finite dimension $d$.
A key direction for future research is to determine if universality-based sharp concentration inequalities can be derived under weaker ergodicity assumptions.

%% file: application.tex
In this section, we apply our results to covariance estimation and principal component analysis (PCA) problems
on Markovian data, to showcase their effectiveness.

We use the same notation as in Section~\ref{subsection:markov_matrix_rosenthal}.
Suppose $(Z_k)_{k\in\NB}$ is a uniformly geometrically ergodic stationary Markov chain (Assumption~\ref{assumption:uniform_geo_ego}) and $\bm{f}\colon\gZ\to\RB^d$ is a bounded vector-valued function.
We first consider the problem of covariance estimation on Markovian data, which is also considered in \citet{fan2021hoeffding}. 
That is, we are interested in analyzing the spectral norm error in estimating the covariance matrix $\bSigma=\mu(\bm{f}\bm{f}^\top)=\EB_\mu[\bm{f}(Z_0)\bm{f}(Z_0)^{\top}]$ using the sample covariance matrix $\hat{\bSigma}_n=n^{-1}\sum_{k=0}^{n-1}\bm{f}(Z_k)\bm{f}(Z_k)^{\top}$.
Let $\bF(z)\coloneq \bm{f}(z)\bm{f}(z)^{\top}-\bSigma$. Then $\bF$ satisfies $\mu(\bF)=\bm{0}_{d{\times}d}$ and $\|\|\bF\|\|_{\infty}\leq \|\|\bm{f}\|\|_{\infty}^2$.
Using our Bernstein inequality for Markov chains (Corollary~\ref{cor:markov_bernstein_uniform}), we obtain that for any $\delta\in(0,1)$, with probability at least $1-\delta$, 
\begin{equation}\label{eq:cov_error}
    \frac{1}{n}\sup_{k\leq n}\norm{\sum_{j=0}^{k-1}\prn{\bm{f}(Z_j)\bm{f}(Z_j)^{\top}-\bSigma}}
    \leq D_5\sqrt{\frac{\log(ed/\delta)\norm{\bSigma_\mu(\bF)}}{n}}
    +D_6\frac{\log(ed/\delta)\tmix\|\|\bm{f}\|\|_\infty^2}{n},
\end{equation}
where $\bSigma_\mu(\bF)$ is the long-run variance of $\bF$ defined in Eqn.~\eqref{eq:long_term_variance}.

Next, we consider the problem of principal component analysis on Markovian data, which is also considered in \citet{neeman2024concentration,kumar2023streaming}.
We denote by $\lambda_1\geq \lambda_2\geq\ldots\geq \lambda_d\geq 0$ the eigenvalues of $\bSigma$, and $\bv_1,\bv_2,\ldots,\bv_d$ the corresponding eigenvectors.
We assume $\lambda_1-\lambda_2>0$.
In the PCA problem, we are interested in analyzing the $\sin^2$ error in estimating the leading eigenvector $\bv_1$ of $\bSigma$ using the leading eigenvector $\hat{\bv}_{1,n}$ of $\hat{\bSigma}_n$. 
By Eqn.~\eqref{eq:cov_error} and Wedin's Theorem \citep{wedin1972perturbation}, we immediately have
\begin{equation*}
\begin{aligned}
    1-\prn{\bv_1^\top \hat{\bv}_{1,n}}^2
    \leq C\bigg\{&\frac{\log(ed/\delta)\norm{\bSigma_\mu(\bF)}}{(\lambda_1-\lambda_2)^2n}+\frac{\log^2(ed/\delta)\tmix^2\|\|\bm{f}\|\|_\infty^4}{(\lambda_1-\lambda_2)^2n^2}\bigg\},
\end{aligned}
\end{equation*}
where $C>0$ is a universal constant.

Similar results hold for the more general geometrically $V$-ergodic Markov chains and unbounded functions $\bm{f}$, and we omit them for brevity.

Compared with the results in \citet{kumar2023streaming,neeman2024concentration}, the leading term in our bounds matches the Markov chain CLT, while the residual term has the sharp $\tmix/n$ dependence. 
Moreover, the Poisson equation result extends to the more general geometric $V$-ergodicity assumption and unbounded functions.

%% file: proof_outline_burkholder.tex
In this section, we outline the proofs of our main results.

\subsection{Proof of the Inequalities for Matrix Martingales}\label{Subsection:proof_outline_burkholder}

In this subsection, our aim is to prove our Rosenthal-Burkholder inequality for the spectral norm of matrix martingales (Theorem~\ref{thm:martingale_rosenthal}).

Based on the Bennett inequality for matrix martingales (Lemma~\ref{lem:matrix_bennet}), we can obtain the following good $\lambda$ inequality for conditionally symmetric matrix martingales, which is a matrix version of \citet[Lemma~4.2]{pinelis1994optimum}.
Before that, we need to introduce the concept of conditionally symmetric martingales.
\begin{definition}[Conditionally Symmetric Martingales]
The martingale $(\bM_k)_{k\in\NB}$ is said to be conditionally symmetric if the conditional distributions of $\bX_k$ and $-\bX_k$ given $\sF_{k-1}$ coincide almost surely for any $k\in\NB$.
\end{definition}
With the assumption that the martingale $(\bM_k)_{k\in\NB}$ is conditionally symmetric, the analysis can be greatly simplified, because the truncated martingale remains a martingale. 
Later on, we will introduce how to perform symmetrization to extend the conclusion to general martingales.
\begin{theorem}[Good $\lambda$ Inequality]\label{thm:good_lambda_inequality}
    Suppose that $(\bM_k)_{k\in\NB}$ is conditionally symmetric.
    For any $\lambda,\delta_1,\delta_2>0$ and $\beta>1+\delta_2$, letting $N=[\beta-(1+\delta_2)]/{\delta_2}$, we have that
    \begin{equation*}
        \PB\prn{\sup_{k\in\NB}\|\bM_k\|>\beta\lambda,\  \norm{\qv{\bM}_\infty}^{1/2}\leq \delta_1\lambda,\ \sup_{k\in\NB}\|\bX_k\|\leq \delta_2\lambda }\leq 2d \prn{\frac{e}{N}\frac{\delta_1^2}{\delta_2^2}}^N \PB\prn{\sup_{k\in\NB}\|\bM_k\|>\lambda}.
    \end{equation*}
    Furthermore, if there exists a $\bUp\in\HB^d_+$ such that $\langle\bM\rangle_\infty\preccurlyeq \bUp$ almost surely, then the above inequality also holds in which $d$ is replaced by $e\cdot r(\bUp\wedge(\delta_1^2\lambda^2))$.
\end{theorem}
The proof is similar to \citet[Lemma~4.2]{pinelis1994optimum} and can be found in Appendix~\ref{appendix:proof_good_lambda_inequality}.
The key ingredient is the Bennett inequality for matrix martingales (Lemma~\ref{lem:matrix_bennet}).
Roughly speaking, this good $\lambda$ inequality indicates that if the supremum of a conditionally symmetric martingale is sufficiently large, then with high probability, both the quadratic variation and the supremum of the martingale difference will also be large.

We observe that the dimension-free factor $r\prn{\bUp\wedge(\delta_1^2\lambda^2)}$ in Theorem~\ref{thm:good_lambda_inequality} increases to $d$ as $\lambda$ decreases to $0$. 
Therefore, to obtain a valid dimension-free result, it is necessary to truncate the error level $\lambda$ before applying the good $\lambda$ inequality.

The following lemma, the truncated $\Phi$ inequality, translates a truncated good $\lambda$ inequality into a moment inequality, which is an extension of \citet[Lemma~7.1]{burkholder1973distribution}.
\begin{lemma}[Truncated $\Phi$ Inequality]\label{lem:truncated_phi_ineq}
Suppose that
$\Phi\colon \RB_+\to\RB_+$ is non-decreasing continuous satisfying $\Phi(\beta\lambda)\leq \gamma \Phi(\lambda)$ for all $\lambda\geq 0$, where $\beta>1$ and $\gamma>0$ are some constants; and $Y, Z$ are nonnegative random variables satisfying 
the following good $\lambda$ inequality
\begin{equation*}
    \PB\prn{Y>\beta \lambda,Z\leq\lambda}\leq\epsilon\PB\prn{Y>\lambda},\quad\forall\lambda\geq \frac{a}{\beta}
\end{equation*}
for some constants $a\geq 0$ and $\epsilon\in(0, 1/\gamma)$.
Then we have that
\begin{equation*}
     \EB\brk{\Phi(Y)}\leq\frac{1}{1-\gamma\epsilon}\Phi(a)+\frac{\gamma}{1-\gamma\epsilon}\EB\brk{\Phi(Z)}.
\end{equation*}
\end{lemma}
The proof can be found in Appendix~\ref{appendix:proof_truncated_phi_inequality}.

Now, we are ready to give a full proof of the Rosenthal-Burkholder inequality (Theorem~\ref{thm:martingale_rosenthal}).
As in the proof of the scalar Rosenthal-Burkholder inequality \citep{hitczenko1990best,pinelis1994optimum}, we first deal with conditionally symmetric martingales, and then generalize the results to general martingales.

\paragraph{Step 1: Conditionally Symmetric Martingales}
Assume that $(\bM_k)_{k\in\NB}$ is conditionally symmetric.
We start with the dimension-dependent case.
Applying the good $\lambda$ inequality (Theorem~\ref{thm:good_lambda_inequality}) and the truncated $\Phi$ inequality (Lemma~\ref{lem:truncated_phi_ineq}) with $Y=\sup_{k\in\NB}\|\bM_k\|$, $Z=\delta_1^{-1}\norm{\qv{\bM}_\infty}^{1/2}+\delta_2^{-1}\sup_{k\in\NB}\|\bX_k\|$, $\Phi(x)=x^p$, $\beta>1$ to be determined, $\gamma=\beta^p$, $a=0$, and 
\begin{equation*}
    \epsilon=2d \prn{\frac{e}{N}\frac{\delta_1^2}{\delta_2^2}}^N,
\end{equation*}
we have the following inequality
    \begin{equation*}
    \norm{\sup_{k\in\NB}\|\bM_k\|}_p\leq \prn{\frac{1}{1-\gamma\epsilon}}^{1/p}\beta \prn{\delta_1^{-1}\norm{\norm{\qv{\bM}_\infty}^{1/2}}_p+\delta_2^{-1}\norm{\sup_{k\in\NB}\norm{\bX_k}}_p },
\end{equation*}
provided that
\begin{equation}\label{eq:rosenthal_parameter_ineq}
    \gamma \epsilon= 2d \prn{\frac{e}{N}\frac{\delta_1^2}{\delta_2^2}}^N\beta^p<1,\quad\text{ where }\ N=\frac{\beta-(1+\delta_2)}{\delta_2}.
\end{equation}
We can choose the following parameters such that $\beta>1+\delta_2$ and the inequality~\eqref{eq:rosenthal_parameter_ineq} holds:
\begin{equation*}
    \beta=1+e^{-p/c}+\frac{1}{c\vee\log d},\quad\delta_1=\frac{1}{5\sqrt{c\vee\log d} \ e^{p/c}},\quad \delta_2=\frac{1}{5(c\vee \log d)},
\end{equation*}
where $c\in [1, p]$ is a parameter to be determined.
Now, we check that inequality~\eqref{eq:rosenthal_parameter_ineq} indeed holds.
First note that $e^{-p/c}\geq e^{-1}$ and then
\begin{equation*}
    \begin{aligned}
        &N=\frac{\beta-1-\delta_2}{\delta_2}=4+5(c\vee \log d)e^{-p/c}\geq 4+5e^{-1}\log d,\\
        &\beta^p=\brk{1+(N+1)\delta_2}^p\leq e^{2pN\delta_2}=e^{2pN/[5(c\vee\log d)]}\leq e^{2pN/(5c)},\\
        &\frac{e}{N}\frac{\delta_1^2}{\delta_2^2}= \frac{e}{4e^{2p/c}(c\vee \log d)^{-1}+5 e^{p/c}}\leq \frac{e}{5}e^{-p/c},
    \end{aligned}
\end{equation*}
hence 
\begin{equation*}
    \begin{aligned}
        \gamma\epsilon=2d \prn{\frac{e}{N}\frac{\delta_1^2}{\delta_2^2}}^N\beta^p&\leq 2d \prn{\frac{e}{5}}^{4+5e^{-1}\log d}e^{-21p/(5c)}< 2(e/5)^4 e^{-21/5}<1,
    \end{aligned}
\end{equation*}
which is the desired conclusion.
Here we used the fact $(e/5)^{5/e}<e^{-1}$
and $2(e/5)^4<1/2$.
Then we have 
\begin{equation*}
   \prn{\frac{1}{1-\gamma\epsilon}}^{1/p}<\frac{1}{\sqrt{1-2(e/5)^4 e^{-21/5}}}  < 1.002.
\end{equation*}

Up to this point, we have shown that for the conditionally symmetric martingale $\{\bM_k\}_{k\in\NB}$, the following inequality holds:
\begin{equation}\label{eq:rosenthal_burkholder_symmetric}
    \norm{\sup_{k\in\NB}\|\bM_k\|}_p\leq 5.01\!\prn{1{+}e^{-p/c}{+}\frac{1}{c{\vee}\log d}}\prn{\sqrt{c{\vee}\log d}\ e^{p{/}c} \norm{\norm{\qv{\bM}_\infty}^{1{/}2}}_p{+}(c{\vee} \log d)\norm{\sup_{k\in\NB}\norm{\bX_k}}_p }.
\end{equation}
Setting $c = p$, we obtain the inequality in the form given in Theorem~\ref{thm:martingale_rosenthal}.
Alternatively, we can choose $c$ as the unique solution to equation $\sqrt{c} = e^{p/c}$, which yields another common form of the Rosenthal-Burkholder inequality $\|\sup_{k\in\NB}\|\bM_k\|\|_p\lesssim [(p/\log p)\vee \log d]( \|\norm{\qv{\bM}_\infty}^{1/2} \|_p  + \| \sup_{k\in\NB}\|\bX_k\|\|_p)$. 
We shall not delve into this form in further detail in this paper.

\textit{Dimension-Free: }To extend the result to a dimension-free version, we only need to replace $a=0$ with
\begin{equation*}
    \tilde{a}=\beta \prn{\delta_1^{-1}\norm{\norm{\qv{\bM}_\infty}^{1/2}}_p+\delta_2^{-1}\norm{\sup_{k\in\NB}\norm{\bX_k}}_p },
\end{equation*}
and $d$ with 
\begin{equation*}
    \tilde{d} =e\cdot r\prn{\bUp\wedge \brk{\norm{\norm{\qv{\bM}_\infty}}_{\frac{p}{2}}}}.
\end{equation*}
And we need to double the right-hand side of the inequality~\eqref{eq:rosenthal_burkholder_symmetric} due to the additional appearance of $\Phi(a)$, that is, if $\langle\bM\rangle_\infty\preccurlyeq \bUp$ almost surely, then it holds that
\begin{equation*}
    \norm{\sup_{k\in\NB}\|\bM_k\|}_p\leq 10.02\!\prn{1{+}e^{-p/c}{+}\frac{1}{c{\vee}\log 
    \tilde d}}\prn{\sqrt{c{\vee}\log \tilde d}\ e^{p{/}c} \norm{\norm{\qv{\bM}_\infty}^{1{/}2}}_p{+}(c{\vee} \log \tilde d)\norm{\sup_{k\in\NB}\norm{\bX_k}}_p }.
\end{equation*}

\paragraph{Step 2: Generalize to Arbitrary Martingales}
In order to generalize the result under the conditionally symmetry assumption to arbitrary martingales, in the special case of independence increment, it suffices to introduce i.i.d. Rademacher random variables; while in the case of martingales, a more sophisticated technique, named decoupling tangent sequence \citep{kwapien1989tangent,hitczenko1988comparison,de2012decoupling}, is required.
The proof in this part is similar to those of the scalar Rosenthal-Burkholder inequality \citep{hitczenko1990best}, and we defer the presentation to Appendix~\ref{appendix:omitted_proof_burkholder}.

%% file: proof_outline_markov.tex
\subsection{Proof of the Inequalities for Markov Chains}\label{Subsection:proof_outline_markov}

In this subsection, we prove the Rosenthal-type inequalities under the spectral norm for Markov chains (Theorem~\ref{thm:markov_rosenthal} and Theorem~\ref{thm:markov_rosenthal_geo_V}).

We consider the Poisson equation with the matrix forcing function $\bF$:
\begin{equation}\label{eq:Poisson_equation}
    \bG-(\rQ\bG)=\bF.
\end{equation}
The following lemma provides the solution and the estimates that will be used below. Its proof can be found in Appendix~\ref{appendix:proof_poisson_equation}.

\begin{lemma}[Solution to the Poisson Equation]\label{lem:poisson_equation}
Consider the matrix Poisson equation~\eqref{eq:Poisson_equation}. Then
\begin{enumerate}[(a)]
    \item Under Assumption~\ref{assumption:uniform_geo_ego} and $\|\|\bF\|\|_{\infty}<\infty$, $\bG=\sum_{k=0}^\infty\rQ^k\bF$ is a solution to the Poisson equation and satisfies
    \begin{equation*}
        \|\|\bG\|\|_{\infty}\leq \frac{8}{3}\tmix\|\|\bF\|\|_{\infty}.
    \end{equation*}
    Furthermore, if there exists an $\bUp\in\HB^d_+$ such that $\bF^2(z)\preccurlyeq\bUp$ for all $z\in\gZ$, then $\bG^2(z)\preccurlyeq(64/9)\tmix^2\bUp$ for all $z\in\gZ$.
    \item Under Assumption~\ref{assumption:V_ego} and $\|\|\bF\|\|_{V^\alpha}<\infty$ for some $\alpha\in(0,1]$, $\bG=\sum_{k=0}^\infty\rQ^k\bF$ is a solution to the Poisson equation and satisfies
    \begin{equation*}
        \norm{\bG(z)}\leq \frac{\varkappa^\alpha2^{3-\alpha}}{3\alpha}\prn{\mu(V)+V(z)}^\alpha\tmix\|\|\bF\|\|_{V^\alpha},\qquad z\in\gZ.
    \end{equation*}
\end{enumerate}
\end{lemma}

We first state a moment inequality variant of the matrix Bernstein inequality in \citet{neeman2024concentration}.

\begin{lemma}[Bernstein Moment Inequality for Markov Chains with a Spectral Gap]\label{lem:constant_gap_bernstein}
Let $P$ be a Markov kernel on $(\gZ,\sZ)$ with stationary distribution $\mu$, and let $\Pi f\coloneq \mu(f)\mathbf{1}$. Suppose that
\begin{equation*}
    \norm{P-\Pi}_{L^2(\mu)\to L^2(\mu)}\leq\frac{1}{2},
\end{equation*}
and $(Y_k)_{k\in\NB}$ is a stationary Markov chain with transition kernel $P$. If a measurable matrix function $\bPhi\colon\gZ\to\HB^d$ satisfies
\begin{equation*}
    \mu(\bPhi)=\bm{0}_{d\times d},\qquad
    \|\|\bPhi\|\|_\infty\leq U,\qquad
    \norm{\mu(\bPhi^2)}\leq\sigma^2,
\end{equation*}
then, for any $r\geq1$ and $N\geq1$,
\begin{equation}\label{eq:constant_gap_bernstein_moment}
    \norm{\norm{\sum_{k=0}^{N-1}\bPhi(Y_k)}}_{\mu,r}
    \leq C\brc{\sqrt{N(r\vee\log d)\sigma^2}+(r\vee\log d)U},
\end{equation}
where $C>0$ is a universal constant.
\end{lemma}

\begin{proof}
By \citet[Theorem~2.6 and Corollary~2.7]{neeman2024concentration}, applied to $\bPhi$ and $-\bPhi$, there exists a universal constant $c>0$ such that
\begin{equation}\label{eq:constant_gap_bernstein_tail}
    \PB_\mu\prn{\norm{\sum_{k=0}^{N-1}\bPhi(Y_k)}\geq u}
    \leq4d^2\exp\prn{-\frac{cu^2}{N\sigma^2+Uu}},\qquad u>0.
\end{equation}
Indeed, the factors $(1+\lambda)/(1-\lambda)$ and $(1-\lambda)^{-1}$ in their result are bounded by universal constants when $\lambda=\norm{P-\Pi}_{L^2(\mu)\to L^2(\mu)}\leq1/2$.

Let $a_d=\log(4d^2)$. It follows from Eqn.~\eqref{eq:constant_gap_bernstein_tail} that, for every $s\geq0$,
\begin{equation*}
    \PB_\mu\brk{\norm{\sum_{k=0}^{N-1}\bPhi(Y_k)}
    >C\brc{\sqrt{N\sigma^2(a_d+s)}+U(a_d+s)}}\leq e^{-s}.
\end{equation*}
Integrating the tail and using $\Gamma(r+1)^{1/r}\leq r$ and $\Gamma(r/2+1)^{1/r}\leq\sqrt r$, we obtain
\begin{equation*}
    \norm{\norm{\sum_{k=0}^{N-1}\bPhi(Y_k)}}_{\mu,r}
    \leq C\brc{\sqrt{N\sigma^2(a_d+r)}+U(a_d+r)}.
\end{equation*}
Since $a_d+r\lesssim r\vee\log d$ for $r\geq1$, Eqn.~\eqref{eq:constant_gap_bernstein_moment} follows.
\end{proof}

Uniform geometric ergodicity gives the required spectral gap after passing to a mixing-time skeleton chain.

\begin{lemma}[Spectral Gap of the Mixing-Time Skeleton Chains]\label{lem:skeleton_gap}
Under Assumption~\ref{assumption:uniform_geo_ego},
\begin{equation}\label{eq:skeleton_gap}
    \norm{\rQ^{2\tmix}-\Pi}_{L^2(\mu)\to L^2(\mu)}\leq\frac{1}{2}.
\end{equation}
\end{lemma}

\begin{proof}
By Assumption~\ref{assumption:uniform_geo_ego},
\begin{equation*}
    \sup_{z\in\gZ}\norm{\rQ^{2\tmix}(z,\cdot)-\mu}_{\TV}
    \leq2\prn{\frac14}^2=\frac18.
\end{equation*}
Thus,
\begin{equation*}
    \norm{\rQ^{2\tmix}-\Pi}_{L^\infty(\mu)\to L^\infty(\mu)}\leq\frac18.
\end{equation*}
On the other hand, the invariance of $\mu$ and the Markov property give
\begin{equation*}
    \norm{\rQ^{2\tmix}f}_{L^1(\mu)}\leq\mu(\rQ^{2\tmix}|f|)=\norm{f}_{L^1(\mu)},
    \qquad
    \norm{\Pi f}_{L^1(\mu)}\leq\norm{f}_{L^1(\mu)}.
\end{equation*}
Hence $\norm{\rQ^{2\tmix}-\Pi}_{L^1(\mu)\to L^1(\mu)}\leq2$. The Riesz--Thorin interpolation theorem now yields Eqn.~\eqref{eq:skeleton_gap}.
\end{proof}

We use the Poisson equation to define
\begin{equation}\label{eq:markov_martingale_decomposition}
    \bX_k\coloneq \bG(Z_k)-(\rQ\bG)(Z_{k-1}),\qquad
    \bM_k\coloneq \sum_{j=1}^k\bX_j,
\end{equation}
so that $\bS_k=\bM_k+\bG(Z_0)-\bG(Z_k)$.
Then $(\bM_k)_{k\in\NB}$ is an $(\sF_k)_{k\in\NB}$-adapted matrix martingale. We define
\begin{equation}\label{eq:def_H_markov}
    \bH\coloneq \rQ\bG^2-(\rQ\bG)^2.
\end{equation}
By the Markov property, $\EB_{k-1}[\bX_k^2]=\bH(Z_{k-1})$.
Moreover, the invariance of $\mu$ and Eqn.~\eqref{eq:Poisson_equation} imply
\begin{equation}\label{eq:H_mean_long_run}
\begin{aligned}
    \mu(\bH)
    &=\mu\prn{\bG^2-(\bG-\bF)^2}\\
    &=\mu\prn{\bF\bG+\bG\bF-\bF^2}\\
    &=\mu\prn{\bF^2}+\sum_{k=1}^\infty\mu\prn{\bF\rQ^k\bF+(\rQ^k\bF)\bF}
    =\bSigma_\mu(\bF).
\end{aligned}
\end{equation}

Now, we are ready to prove Theorem~\ref{thm:markov_rosenthal}.
\begin{proof}[Proof of Theorem~\ref{thm:markov_rosenthal}]
We first assume that $Z_0\sim\mu$. Consider the skeleton chain and the martingale sampled at its block endpoints,
\begin{equation*}
    Y_j\coloneq Z_{2j\tmix},\qquad
    \wtilde{\bM}_j\coloneq \bM_{2j\tmix},\qquad j\in\NB.
\end{equation*}
The chain $(Y_j)_{j\in\NB}$ is stationary with transition kernel $\rQ^{2\tmix}$, and $(\wtilde{\bM}_j)_{j\in\NB}$ is adapted to $(\sF_{2j\tmix})_{j\in\NB}$. Its martingale differences satisfy
\begin{equation}\label{eq:block_martingale_difference}
\begin{aligned}
    \wtilde{\bX}_j
    &\coloneq \wtilde{\bM}_j-\wtilde{\bM}_{j-1}\\
    &=\sum_{r=2(j-1)\tmix}^{2j\tmix-1}\bF(Z_r)
      +\bG(Z_{2j\tmix})-\bG(Z_{2(j-1)\tmix}),\qquad j\geq1.
\end{aligned}
\end{equation}
By Lemma~\ref{lem:poisson_equation} (a),
\begin{equation}\label{eq:block_increment_bound}
    \norm{\wtilde{\bX}_j}
    \leq2\tmix\norm{\norm{\bF}}_\infty+2\norm{\norm{\bG}}_\infty
    \leq\frac{22}{3}\tmix\norm{\norm{\bF}}_\infty.
\end{equation}

For $z\in\gZ$, define
\begin{equation}\label{eq:def_block_K}
    \bK(z)\coloneq \EB_z[\wtilde{\bX}_1^2].
\end{equation}
The Markov property gives
\begin{equation}\label{eq:block_quadratic_variation}
    \qv{\wtilde{\bM}}_N=\sum_{j=0}^{N-1}\bK(Y_j).
\end{equation}
To identify $\bK$, note that $\wtilde{\bX}_1=\sum_{\ell=1}^{2\tmix}\bX_\ell$. For $1\leq \ell<r\leq2\tmix$,
\begin{equation*}
    \EB_z[\bX_\ell\bX_r]
    =\EB_z[\bX_\ell\EB_{r-1}(\bX_r)]=\bm{0}_{d\times d},
\end{equation*}
and similarly $\EB_z[\bX_r\bX_\ell]=\bm{0}_{d\times d}$. Therefore, the cross terms vanish without any commutativity assumption, and
\begin{equation}\label{eq:block_K_representation}
    \bK(z)=\sum_{\ell=1}^{2\tmix}\EB_z[\bX_\ell^2]
    =\sum_{\ell=0}^{2\tmix-1}(\rQ^\ell\bH)(z).
\end{equation}
It follows from Eqn.~\eqref{eq:H_mean_long_run} that
\begin{equation}\label{eq:block_K_mean}
    \mu(\bK)=2\tmix\bSigma_\mu(\bF).
\end{equation}

Equation~\eqref{eq:block_increment_bound} also gives
\begin{equation}\label{eq:block_K_envelope}
    \bm{0}_{d\times d}\preccurlyeq\bK(z)
    \preccurlyeq\prn{\frac{22}{3}}^2\tmix^2\norm{\norm{\bF}}_\infty^2\bI_d.
\end{equation}
Let $\bar\bK\coloneq \bK-\mu(\bK)$. Then $\mu(\bar\bK)=\bm{0}_{d\times d}$ and
\begin{equation}\label{eq:bar_K_envelope}
    \|\|\bar\bK\|\|_\infty
    \leq\prn{\frac{22}{3}}^2\tmix^2\norm{\norm{\bF}}_\infty^2.
\end{equation}
And $\mu(\bar\bK^2)$ can be upper-bounded by
\begin{equation}\label{eq:bar_K_second_moment}
\begin{aligned}
    \mu(\bar\bK^2)
    &=\mu(\bK^2)-\mu(\bK)^2\\
    &\preccurlyeq\mu(\bK^2)\\
    &\preccurlyeq\prn{\frac{22}{3}}^2\tmix^2\norm{\norm{\bF}}_\infty^2\mu(\bK)\\
    &=2\prn{\frac{22}{3}}^2\tmix^3\norm{\norm{\bF}}_\infty^2\bSigma_\mu(\bF).
\end{aligned}
\end{equation}
Here the second inequality follows from $\bK^2\preccurlyeq(22/3)^2\tmix^2\|\|\bF\|\|_\infty^2\bK$, which is a consequence of Eqn.~\eqref{eq:block_K_envelope}.

By Lemma~\ref{lem:skeleton_gap}, we may apply Lemma~\ref{lem:constant_gap_bernstein} to $(Y_j)_{j\in\NB}$ and $\bar\bK$. Taking $r=p/2$ in Eqn.~\eqref{eq:constant_gap_bernstein_moment} and using Eqns.~\eqref{eq:bar_K_envelope} and \eqref{eq:bar_K_second_moment}, we obtain
\begin{equation}\label{eq:block_K_concentration}
\begin{aligned}
    \norm{\norm{\sum_{j=0}^{N-1}\bar\bK(Y_j)}}_{\mu,p/2}
    \leq C\bigg\{&\tmix\norm{\norm{\bF}}_\infty
    \sqrt{N\tmix(p\vee\log d)\norm{\bSigma_\mu(\bF)}}+(p\vee\log d)\tmix^2\norm{\norm{\bF}}_\infty^2\bigg\}.
\end{aligned}
\end{equation}
Combining Eqns.~\eqref{eq:block_quadratic_variation}, \eqref{eq:block_K_mean}, and \eqref{eq:block_K_concentration} yields
\begin{equation}\label{eq:block_quadratic_complete_square}
\begin{aligned}
    \norm{\norm{\qv{\wtilde{\bM}}_N}}_{\mu,p/2}
    &\leq2N\tmix\norm{\bSigma_\mu(\bF)}
      +\norm{\norm{\sum_{j=0}^{N-1}\bar\bK(Y_j)}}_{\mu,p/2}\\
    &\leq C\brc{\sqrt{N\tmix\norm{\bSigma_\mu(\bF)}}
      +\tmix\sqrt{p\vee\log d}\norm{\norm{\bF}}_\infty}^2.
\end{aligned}
\end{equation}
Consequently,
\begin{equation}\label{eq:block_quadratic_root_bound}
    \norm{\norm{\qv{\wtilde{\bM}}_N}^{1/2}}_{\mu,p}
    \leq C\brc{\sqrt{N\tmix\norm{\bSigma_\mu(\bF)}}
      +\tmix\sqrt{p\vee\log d}\norm{\norm{\bF}}_\infty}.
\end{equation}

We apply Theorem~\ref{thm:martingale_rosenthal} to the stopped martingale $(\wtilde{\bM}_{j\wedge N})_{j\in\NB}$. By Eqns.~\eqref{eq:block_increment_bound} and \eqref{eq:block_quadratic_root_bound},
\begin{equation}\label{eq:block_martingale_bound}
    \norm{\sup_{j\leq N}\norm{\wtilde{\bM}_j}}_{\mu,p}
    \leq C\brc{\sqrt{N\tmix(p\vee\log d)\norm{\bSigma_\mu(\bF)}}
    +(p\vee\log d)\tmix\norm{\norm{\bF}}_\infty}.
\end{equation}

It remains to pass from the block endpoints to all partial sums. Let $N=\lfloor n/(2\tmix)\rfloor$. For every $k\leq n$, set $j=\lfloor k/(2\tmix)\rfloor$. Then
\begin{equation*}
    \norm{\bS_k-\bS_{2j\tmix}}\leq2\tmix\norm{\norm{\bF}}_\infty,
\end{equation*}
while Eqn.~\eqref{eq:markov_martingale_decomposition} gives
\begin{equation*}
    \bS_{2j\tmix}=\wtilde{\bM}_j+\bG(Z_0)-\bG(Z_{2j\tmix}).
\end{equation*}
Hence,
\begin{equation}\label{eq:block_to_all_times}
    \sup_{k\leq n}\norm{\bS_k}
    \leq\sup_{j\leq N}\norm{\wtilde{\bM}_j}
    +\frac{22}{3}\tmix\norm{\norm{\bF}}_\infty.
\end{equation}
If $n\geq2\tmix$, Eqns.~\eqref{eq:block_martingale_bound} and \eqref{eq:block_to_all_times}, together with $2N\tmix\leq n$, imply
\begin{equation}\label{eq:sharp_stationary_unnormalized}
    \norm{\sup_{k\leq n}\norm{\bS_k}}_{\mu,p}
    \leq C\brc{\sqrt{n(p\vee\log d)\norm{\bSigma_\mu(\bF)}}
    +(p\vee\log d)\tmix\norm{\norm{\bF}}_\infty}.
\end{equation}
If $n<2\tmix$, the same bound follows directly from
\begin{equation*}
    \sup_{k\leq n}\norm{\bS_k}\leq n\norm{\norm{\bF}}_\infty
    \leq2\tmix\norm{\norm{\bF}}_\infty.
\end{equation*}
Dividing Eqn.~\eqref{eq:sharp_stationary_unnormalized} by $n$ proves Eqn.~\eqref{eq:sharp_stationary_rosenthal}.

We next consider an arbitrary initial distribution $\xi\in\gP$. Let $(\wtilde Z_k)_{k\in\NB}$ be a stationary copy of the chain. By Assumption~\ref{assumption:uniform_geo_ego},
\begin{equation*}
    \sup_{z,z'\in\gZ}\norm{\rQ^{\tmix}(z,\cdot)-\rQ^{\tmix}(z',\cdot)}_{\TV}\leq1.
\end{equation*}
Recall that the total variation norm used in this paper takes values in $[0,2]$ for the difference of two probability measures. A maximal coupling at times $\tmix,2\tmix,\ldots$, followed by the same transitions after the two chains meet, therefore gives a coupling time $\tau$ satisfying
\begin{equation}\label{eq:coupling_geometric_tail}
    \PB(\tau>j\tmix)\leq2^{-j},\qquad j\in\NB.
\end{equation}
This is the standard maximal exact coupling construction; see \citet[Lemmas~19.3.6 and~19.3.9]{douc2018markov}. It follows that
\begin{equation*}
    \PB(\tau>s)\leq2\exp\prn{-\frac{(\log2)s}{\tmix}},\qquad s\geq0,
\end{equation*}
and therefore, for $p\geq2$,
\begin{equation}\label{eq:coupling_moment}
    \norm{\tau}_p
    \leq\brk{2\Gamma(p+1)}^{1/p}\frac{\tmix}{\log2}
    \leq\frac{2}{\log2}p\tmix.
\end{equation}
Let $\wtilde{\bS}_k=\sum_{j=0}^{k-1}\bF(\wtilde Z_j)$. Since the chains agree after time $\tau$,
\begin{equation*}
    \sup_{k\leq n}\norm{\bS_k-\wtilde{\bS}_k}
    \leq2\norm{\norm{\bF}}_\infty(\tau\wedge n).
\end{equation*}
Minkowski's inequality and Eqn.~\eqref{eq:coupling_moment} give
\begin{equation*}
    \norm{\sup_{k\leq n}\norm{\bS_k}}_{\xi,p}
    \leq\norm{\sup_{k\leq n}\norm{\bS_k}}_{\mu,p}
    +\frac{4}{\log2}p\tmix\norm{\norm{\bF}}_\infty,
\end{equation*}
which proves Eqns.~\eqref{eq:from_non_stationary_to_stationary} and \eqref{eq:sharp_general_rosenthal}.

For the rectangular case, consider the Hermitian dilation
\begin{equation*}
    \sH(\bR)\coloneq \begin{bmatrix}\bm{0}&\bR\\\bR^*&\bm{0}\end{bmatrix}.
\end{equation*}
Then $\norm{\sH(\bR)}=\norm{\bR}$, $\sH(\rQ^k\bR)=\rQ^k\sH(\bR)$, and the long-run variance of $\sH(\bR)$ is block diagonal, with its two diagonal blocks equal to the two matrices in Eqn.~\eqref{eq:long_term_variance_rect}. Applying the self-adjoint result to $\sH(\bR)$ completes the proof.
\end{proof}

\begin{proof}[Proof of Corollaries~\ref{cor:markov_hoeffding_uniform} and \ref{cor:markov_bernstein_uniform}]
Assumption~\ref{assumption:uniform_geo_ego} and $\mu(\bF)=\bm{0}_{d\times d}$ imply
\begin{equation*}
    \norm{(\rQ^k\bF)(z)}
    \leq2\prn{\frac14}^{\lfloor k/\tmix\rfloor}\norm{\norm{\bF}}_\infty.
\end{equation*}
Consequently,
\begin{equation}\label{eq:long_run_variance_uniform_bound}
\begin{aligned}
    \norm{\bSigma_\mu(\bF)}
    &\leq\norm{\norm{\bF}}_\infty^2
      +4\norm{\norm{\bF}}_\infty^2\sum_{k=1}^\infty\prn{\frac14}^{\lfloor k/\tmix\rfloor}\\
    &\leq\frac{19}{3}\tmix\norm{\norm{\bF}}_\infty^2.
\end{aligned}
\end{equation}
Take $p=2\vee\log(1/\delta)$. By Markov's inequality,
\begin{equation*}
    \PB_\xi\brk{\sup_{k\leq n}\norm{\bS_k}
    >e\norm{\sup_{k\leq n}\norm{\bS_k}}_{\xi,p}}\leq e^{-p}\leq\delta.
\end{equation*}
Since $p\vee\log d\leq2\log(ed/\delta)$, Corollary~\ref{cor:markov_bernstein_uniform} follows from Eqn.~\eqref{eq:sharp_general_rosenthal}. Combining it with Eqn.~\eqref{eq:long_run_variance_uniform_bound} proves Corollary~\ref{cor:markov_hoeffding_uniform}.
\end{proof}

We finally turn to geometrically $V$-ergodic chains. We use the following auxiliary bound, whose proof can be found in Appendix~\ref{appendix:proof_crude_rosenthal_markov}.

\begin{lemma}\label{lem:crude_rosenthal_markov}
\begin{enumerate}[(a)]
    \item Under Assumption~\ref{assumption:uniform_geo_ego} and $\|\|\bF\|\|_\infty<\infty$, for any $p\geq2$,
    \begin{equation*}
    \begin{aligned}
        \norm{\sup_{k\leq n}\norm{\bS_k}}_{\mu,p}
        \leq&\frac{16}{3}C_1\sqrt{(p\vee\log d)n\tmix}\norm{\norm{\bF}}_\infty\\
        &+\prn{\frac{16}{3}C_2(p\vee\log d)+\frac{19}{3}}\tmix\norm{\norm{\bF}}_\infty.
    \end{aligned}
    \end{equation*}
    \item Under Assumption~\ref{assumption:V_ego} and $\|\|\bF\|\|_{V^{1/p}}<\infty$ for $p\geq2$,
    \begin{equation*}
    \begin{aligned}
        \norm{\sup_{k\leq n}\norm{\bS_k}}_{\mu,p}
        \leq&\frac{8}{3}C_1p\sqrt{(p\vee\log d)n\tmix}
        \brk{\varkappa\mu(V)}^{1/p}\norm{\norm{\bF}}_{V^{1/p}}\\
        &+\prn{\frac{16}{3}C_2(p\vee\log d)+\frac{11}{3}}p\tmix n^{1/p}
        \brk{\varkappa\mu(V)}^{1/p}\norm{\norm{\bF}}_{V^{1/p}}.
    \end{aligned}
    \end{equation*}
\end{enumerate}
\end{lemma}

\begin{proof}[Proof of Theorem~\ref{thm:markov_rosenthal_geo_V}]
We use the martingale decomposition~\eqref{eq:markov_martingale_decomposition}. For a stationary chain, Theorem~\ref{thm:martingale_rosenthal}, Eqn.~\eqref{eq:H_mean_long_run}, and Minkowski's inequality give
\begin{equation}\label{eq:V_markov_basic_bound}
\begin{aligned}
    \norm{\sup_{k\leq n}\norm{\bS_k}}_{\mu,p}
    \leq& C_1\sqrt{n(p\vee\log d)}\norm{\bSigma_\mu(\bF)}^{1/2}\\
    &+C_1\sqrt{p\vee\log d}\norm{\norm{\sum_{k=0}^{n-1}\bar\bH(Z_k)}}_{\mu,p/2}^{1/2}\\
    &+C_2(p\vee\log d)\norm{\sup_{k\leq n}\norm{\bX_k}}_{\mu,p}
    +2\norm{\norm{\bG(Z_0)}}_{\mu,p},
\end{aligned}
\end{equation}
where $\bar\bH=\bH-\mu(\bH)$. Lemma~\ref{lem:poisson_equation} (b) implies
\begin{equation*}
    \norm{\norm{\bar\bH}}_{V^{2/p}}^{1/2}
    \lesssim p\brk{\varkappa\mu(V)}^{1/p}\tmix\norm{\norm{\bF}}_{V^{1/p}}.
\end{equation*}
Applying Lemma~\ref{lem:crude_rosenthal_markov} (b) to $\bar\bH$ yields
\begin{equation*}
\begin{aligned}
    &C_1\sqrt{p\vee\log d}\norm{\norm{\sum_{k=0}^{n-1}\bar\bH(Z_k)}}_{\mu,p/2}^{1/2}\\
    \lesssim&\ p^{3/2}(p\vee\log d)^{3/4}\tmix^{5/4}n^{1/4}
    \brk{\varkappa\mu(V)}^{2/p}\norm{\norm{\bF}}_{V^{1/p}}\\
    &+p^{3/2}(p\vee\log d)\tmix^{3/2}n^{1/p}
    \brk{\varkappa\mu(V)}^{2/p}\norm{\norm{\bF}}_{V^{1/p}}.
\end{aligned}
\end{equation*}
The last two terms in Eqn.~\eqref{eq:V_markov_basic_bound} are bounded by the same right-hand side. Dividing by $n$ proves the stationary conclusion.

For an arbitrary $\xi\in\gP_V$, the maximal exact coupling argument in \citet[Section~4]{durmus2023rosenthal} gives
\begin{equation*}
\begin{aligned}
    \frac{1}{n}\norm{\sup_{k\leq n}\norm{\bS_k}}_{\xi,p}
    \leq&\frac{1}{n}\norm{\sup_{k\leq n}\norm{\bS_k}}_{\mu,p}\\
    &+\frac{32}{15n}\varkappa^{1/p}\prn{\xi(V)+\mu(V)}^{1/p}p\tmix
    \norm{\norm{\bF}}_{V^{1/p}},
\end{aligned}
\end{equation*}
which completes the proof.
\end{proof}

%% file: conclusion.tex
In this work, we have studied probability bounds for the spectral norm of sums of dependent random matrices.
First, we have proved a Rosenthal-Burkholder inequality for matrix martingales (Theorem~\ref{thm:martingale_rosenthal}) using the extrapolation method (good $\lambda$ inequality) \citep{hitczenko1990best,hitczenko1994domination,pinelis1994optimum} and the matrix Bennett inequality \citep{tropp2011freedmans, kroshnin2024bernstein}.
By discretization, we have shown a sharp
Burkholder-Davis-Gundy inequality for continuous matrix local martingales (Theorem~\ref{thm:martingale_BDG}) as a corollary of Theorem~\ref{thm:martingale_rosenthal}.
Furthermore, combining this martingale inequality with the Poisson equation approach \citep{durmus2023rosenthal}, a blocking argument, and a matrix Bernstein inequality on the mixing-time skeleton chain \citep{neeman2024concentration}, we have derived matrix Rosenthal inequalities for Markov chains under both the uniform geometric ergodicity assumption and the geometric $V$-ergodicity assumption in Section~\ref{subsection:markov_matrix_rosenthal}.
By the Markov inequality, these moment inequalities can be further translated into concentration inequalities including Hoeffding and Bernstein inequalities.
The sharp uniformly geometrically ergodic bound has the long-run variance in its leading term and the standard $\tmix/n$ dependence in its residual term.

Compared with previous works on matrix concentration inequalities for Markov chains, our results yield a sharp leading term that matches the Markov chain CLT, rather than relying on a sub-optimal variance proxy. 
Our results do not require the assumption of a non-zero absolute spectral gap or $\psi$-mixing condition. 
This expands its applicability to a broader range of problems, such as policy evaluation in reinforcement learning \citep{duan2022policy} and PCA \citep{kumar2023streaming} with Markovian data.
And our results do not need the assumption of the boundedness of functions under the geometric ergodicity assumption. 
The martingale inequalities also have dimension-free versions, enabling applications to infinite-dimensional settings.
The advantages of our results over previous works are demonstrated in the covariance estimation and PCA problems on Markovian data.
And as mentioned in this paper, our results can also be extended to more general geometrically ergodic Markov chains with respect to the Wasserstein semimetric, as well as to geometrically ergodic continuous Markov processes.
One direction for future work is to extend the sharp mixing-time dependence to geometrically $V$-ergodic chains and unbounded matrix functions.

%% file: proof_burkholder.tex
In this section, we complete the proofs omitted in Section~\ref{Subsection:proof_outline_burkholder}.

\subsection{Proof of Theorem~\ref{thm:good_lambda_inequality}}\label{appendix:proof_good_lambda_inequality}
\begin{proof}
We first recall the Bennett inequality for matrix martingales with bounded martingale difference, which is provided in \citet[Theorem~3.1]{tropp2011freedmans}  \citet[Theorem~2.3]{kroshnin2024bernstein}.
\begin{lemma}[Bennett Inequality for Matrix Martingale]\label{lem:matrix_bennet}
Suppose that there exists a $\bUp\in\HB^d_+$ such that $\langle\bM\rangle_\infty\preccurlyeq \bUp$ almost surely, and $\|\bX_k\|\leq U$ for any $k\in\NB$ almost surely for some $U>0$.
Then for any $\varepsilon>0$, it holds that
\begin{equation*}
    \begin{aligned}
       \PB\prn{\sup_{k\in\NB}\norm{\bM_k}\geq \varepsilon}\leq& 2(d\wedge (e\cdot r(\bUp))) \exp\brc{-\frac{\norm{\bUp}}{U^2}h\prn{\frac{U\varepsilon}{\norm{\bUp}}}}\leq 2(d\wedge (e\cdot r(\bUp)))\prn{\frac{e\norm{\bUp}}{\varepsilon U}}^{{\varepsilon}/{U}},
    \end{aligned}
\end{equation*}
where $h(x)=(1+x)\log(1+x)-x$ for any $x>0$.
\end{lemma}

For any $k\in\NB$, let $\bar{\bX}_k=\bX_k\ind\{\|\bX_k\|\leq \delta_2\lambda\}$ be the truncated version of $\bX_k$, $\bar{\bM}_k=\sum_{j=0}^k\bar{\bX}_j$ be the truncated martingale, and $\langle\bar{\bM}\rangle_k=\sum_{j=0}^k \EB_{j-1}\bar{\bX}_j^2 \preccurlyeq \langle\bM\rangle_k$ be the corresponding quadratic variation sequence.
$\bar{\bM}_k$ is indeed a martingale because
\begin{equation*}
    \begin{aligned}
        \EB_{k-1}\brk{\bar{\bX}_k}=&\EB_{k-1}\brk{\bX_k\ind\{\|\bX_k\|\leq \delta_2\lambda\}}\\
        =&\EB_{k-1}\brk{{\bX}_k}-\EB_{k-1}\brk{\bX_k\ind\{\|\bX_k\|> \delta_2\lambda\}}\\
        =&\EB_{k-1}\brk{\bX_k\ind\{\|\bX_k\|> \delta_2\lambda\}}\\
        =&\EB_{k-1}\brk{{\bX}_k}-\EB_{k-1}\brk{\bX_k\ind\{\|\bX_k\|\leq \delta_2\lambda\}}\\
        =&-\EB_{k-1}\brk{\bar{\bX}_k}\\
        =&\bm{0}_{d\times d},
    \end{aligned}
\end{equation*}
where we used the fact that $(\bM_k)_{k\in\NB}$ is a martingale and it is conditionally symmetric.
By definition of the truncated martingale, on the event $\{\sup_{k\in\NB}\|\bX_k\|\leq \delta_2\lambda\}$, we have $\bM_k=\bar{\bM}_k$ for any $k\in\NB$.

Next, we define three stopping times: $\tau_1=\inf\{k\in\NB\colon \|{\bM_k}\|>\lambda \}$, $\tau_2=\inf\{k\in\NB\colon \|{\bM_k}\|>\beta\lambda \}$ ($\tau_1\leq \tau_2$), $\tau_3=\inf\{k\in\NB\colon \|\langle\bar{\bM}\rangle_{k+1}\|^{1/2}>\delta_1\lambda \}$.
By definition of the stopping times, we have
\begin{equation}\label{eq:stopping_time_meaning}
    \brc{\sup_{k\in\NB}\|{\bM}_k\|>\lambda}=\brc{\tau_1<\infty},\quad\brc{\sup_{k\in\NB}\|{\bM}_k\|>\beta\lambda}=\brc{\tau_2<\infty},\quad\brc{\norm{\qv{\bar\bM}_\infty}^{1/2}\leq \delta_1\lambda}=\brc{\tau_3=\infty}.
\end{equation}
We also consider the stopped martingale $\hat{\bM}_k=\bar{\bM}_{(k\wedge\tau_2\wedge\tau_3)\vee \tau_1}-\bar{\bM}_{\tau_1}$ conditionally on the sigma-algebra of the $\tau_1$-past
\begin{equation*}
    \sF_{\tau_1}=\brc{A\in\sF\colon A\cap\{\tau_1\leq k \}\in\sF_k,\quad \forall k\in\NB }.
\end{equation*}
That is $(\hat{\bM}_k)_{k\in\NB}$ is a $(\sF_{k\vee\tau_1})_{k\in\NB}$-adapted martingale
with the martingale difference $\hat{\bX}_k=\hat{\bM}_k-\hat{\bM}_{k-1}=\bar{\bX}_{k}\ind\{\tau_1<k<\tau_2\wedge\tau_3\}$ upper-bounded by $\delta_2\lambda$ and square root of the quadratic variation upper-bounded by $\delta_1\lambda$.

Now, we are ready to derive the good $\lambda$ inequality.
\begin{equation}\label{eq:proof_good_lambda}
    \begin{aligned}
        &\PB\prn{\sup_{k\in\NB}\|\bM_k\|>\beta\lambda,\  \norm{\qv{\bM}_\infty}^{1/2}\leq \delta_1\lambda,\ \sup_{k\in\NB}\|\bX_k\|\leq \delta_2\lambda }\\
        &\qquad\leq\PB\prn{\sup_{k\in\NB}\|{\bM}_k\|>\beta\lambda,\  \norm{\qv{\bar\bM}_\infty}^{1/2}\leq \delta_1\lambda,\sup_{k\in\NB}\|\bX_k\|\leq \delta_2\lambda }\\
        &\qquad\leq\PB\prn{\tau_2<\infty,\tau_3=\infty,\sup_{k\in\NB}\|\bX_k\|\leq \delta_2\lambda  }\\
        &\qquad\leq \PB\prn{\sup_{k\in\NB}\|\hat{\bM}_k\|>(\beta-1-\delta_2)\lambda }\\
        &\qquad= \EB\brk{\PB\prn{ \sup_{k\in\NB}\|\hat{\bM}_k\|>(\beta-1-\delta_2)\lambda\Big|\sF_{\tau_1}} \ind\brc{\tau_1<\infty}}\\
        &\qquad\leq 2d \prn{\frac{e}{N}\frac{\delta_1^2}{\delta_2^2}}^N \EB\brk{\ind\brc{\tau_1<\infty}}\\
        &\qquad= 2d \prn{\frac{e}{N}\frac{\delta_1^2}{\delta_2^2}}^N\PB\prn{\sup_{k\in\NB}\|\bM_k\|>\lambda},
    \end{aligned}
\end{equation}
which is the desired conclusion.
Here, we explain some details in the derivation of Eqn~\eqref{eq:proof_good_lambda}.
In the first inequality, we used the fact that $\langle\bar{\bM}\rangle_\infty\preccurlyeq \langle\bM\rangle_\infty$ almost surely.
In the second inequality, we used Eqn.~\eqref{eq:stopping_time_meaning}.
In the third inequality, we used the fact that on the event $\{\tau_2<\infty,\tau_3=\infty,\sup_{k\in\NB}\|\bX_k\|\leq \delta_2\lambda\}$, it holds that $\sup_{k<\tau_1}\norm{{\bM}_k}=\sup_{k<\tau_1}\norm{\bar{\bM}_k}<\lambda$ and $\norm{\bar{\bM}_{(\tau_2\wedge\tau_3) -1}}>\norm{{\bM}_{\tau_2}}-\norm{\bX_{\tau_2}}>(\beta-\delta_2)\lambda$.
And in the last inequality, we used the matrix Bennett inequality (Lemma~\ref{lem:matrix_bennet}) with error term $\varepsilon=(\beta-1-\delta_2)\lambda$, parameter $U=\delta_2 \lambda$, and 
$\|\langle\hat{\bM}_k\rangle_{\infty} \|\leq\delta_1^2\lambda^2$ given $\gF_{\tau_1}$.

If we further assume that there exists a $\bUp\in\HB^d_+$ such that $\langle\bM\rangle_\infty\preccurlyeq \bUp$ almost surely, then the fourth inequality also holds if we replace $d$ with $e\cdot r(\bUp\wedge(\delta_1^2\lambda^2))$ by the dimension-free matrix Bennett inequality (Lemma~\ref{lem:matrix_bennet}).
\end{proof}

\subsection{Proof of Theorem~\ref{thm:martingale_BDG}}\label{appendix:proof_bdg_inequality}
\begin{proof}
We first consider $\tau=\infty$ and the bounded case: $(\bM_t)_{t\geq 0}$ is a continuous martingale satisfying $\sup_{t\geq 0}\norm{\bM_t}\le K, \norm{\qv{\bM}_\infty}\le K$ for some constant $K>0$. 
For a fixed discretization level $n\in\NB$, we define $\bM_k^{(n)}=\bM_{\frac{k}{2^n}}$ for any $k\in\NB$ ($\bM_{-1}^{(n) }:=\bm{0}_{d{\times}d}$). 
Apply Theorem~\ref{thm:martingale_rosenthal} to $(\bM_k^{(n)})_{k\in\NB}$, we have
\begin{equation}
    \begin{aligned}
       &\norm{\sup_{k\in\NB}\|\bM_{\frac{k}{2^n}}\|}_p\\
       \leq &C_{R,1}\sqrt{p\vee \log d}\norm{\norm{\sum_{k=0}^\infty\EB\brk{(\bM_{\frac{k}{2^n}}-\bM_{\frac{k-1}{2^n}})^2\big|{\sF_{\frac{k-1}{2^n}}}}}^{1/2} }_p + C_{R,2}(p\vee\log d)\norm{ \sup_{k\in\NB}\|\bM_{\frac{k}{2^n}}-\bM_{\frac{k-1}{2^n}}\|}_p.
    \end{aligned}
\end{equation}
Now we take the limit of each part of the above inequality as $n\to \infty$.
By Beppo Levi's Theorem, the left hand side converges to $\norm{\sup_{t\geq 0}\|\bM_t\|}_p$. 
By dominated convergence theorem, the second term of the right hand side converges to $0$. It suffices to prove $\lim_{n\to \infty}\sum_{k=0}^{\infty}\EB\brk{(\bM_{\frac{k}{2^n}}-\bM_{\frac{k-1}{2^n}})^2\big|{\sF_{\frac{k-1}{2^n}}}}=\qv{\bM}_\infty$ in $L^{\frac{p}{2}}$. 
Here the convergence of matrices is defined entry-wise.

\begin{lemma}\label{lem:square_of_variance}
\begin{equation}
    \begin{aligned}       \EB\norm{\sum_{k=0}^{\infty}(\bM_{\frac{k}{2^n}}-\bM_{\frac{k-1}{2^n}})^2}^2\le 48d^4 K^4.
    \end{aligned}
\end{equation}
\end{lemma}
\begin{proof}
    \begin{equation}
    \begin{aligned}       \EB\norm{\sum_{k=0}^{\infty}(\bM_{\frac{k}{2^n}}-\bM_{\frac{k-1}{2^n}})^2}^2
    &\le\EB\prn{\tr\prn{{\sum_{k=0}^{\infty}(\bM_{\frac{k}{2^n}}-\bM_{\frac{k-1}{2^n}})^2}}}^2\\
    &=\EB\prn{\sum_{i,j=1}^{d}{\sum_{k=0}^{\infty}\prn{[\bM_{\frac{k}{2^n}}]_{ij}-[\bM_{\frac{k-1}{2^n}}]_{ij}}^2}}^2\\
    &\le d^2\sum_{i,j=1}^{d}\EB\prn{{\sum_{k=0}^{\infty}\prn{[\bM_{\frac{k}{2^n}}]_{ij}-[\bM_{\frac{k-1}{2^n}}]_{i,j}}^2}}^2.
    \end{aligned}
\end{equation}
By \citep[Lemma~1.5.9]{karatzas1988brownian}, each term in the last summation is upper bounded by $48 K^4$.
\end{proof}

\begin{lemma}\label{lem:fourth_moment}
\begin{equation}
    \begin{aligned}
       \lim_{n\to \infty}\EB\norm{\sum_{k=0}^{\infty}(\bM_{\frac{k}{2^n}}-\bM_{\frac{k-1}{2^n}})^4}=0
    \end{aligned}
\end{equation}
\end{lemma}
\begin{proof}
\begin{equation}
    \begin{aligned}
       \EB\norm{\sum_{k=0}^{\infty}(\bM_{\frac{k}{2^n}}-\bM_{\frac{k-1}{2^n}})^4}
       &\le \EB\norm{\sum_{k=0}^{\infty}\prn{\sup_{k\in\NB}\|\bM_{\frac{k}{2^n}}-\bM_{\frac{k-1}{2^n}}\|^2}(\bM_{\frac{k}{2^n}}-\bM_{\frac{k-1}{2^n}})^2}\\
       &= \EB\prn{\sup_{k\in\NB}\|\bM_{\frac{k}{2^n}}-\bM_{\frac{k-1}{2^n}}\|^2}\norm{\sum_{k=0}^{\infty}(\bM_{\frac{k}{2^n}}-\bM_{\frac{k-1}{2^n}})^2}\\
        &\le \prn{\EB\sup_{k\in\NB}\|\bM_{\frac{k}{2^n}}-\bM_{\frac{k-1}{2^n}}\|^4}^{{1}/{2}}\prn{\EB\norm{\sum_{k=0}^{\infty}(\bM_{\frac{k}{2^n}}-\bM_{\frac{k-1}{2^n}})^2}^2}^{{1}/{2}}\\
        &\le \sqrt{48}d^2K^2\prn{\EB\sup_{k\in\NB}\|\bM_{\frac{k}{2^n}}-\bM_{\frac{k-1}{2^n}}\|^4}^{{1}/{2}}.
    \end{aligned}
\end{equation}
    In the last step, we use the upper bound in Lemma~\ref{lem:square_of_variance}. 
    The expectation in the last step converges to $0$ by dominated convergence theorem.
\end{proof}

\begin{lemma}\label{lem:dis_quadratic_converge}
\begin{equation}
    \begin{aligned}
       \lim_{n\to \infty}\sum_{k=0}^{\infty}\EB\brk{(\bM_{\frac{k}{2^n}}-\bM_{\frac{k-1}{2^n}})^2\big|{\sF_{\frac{k-1}{2^n}}}}=\qv{\bM}_\infty\ {\rm in} \ L^{2}.
    \end{aligned}
\end{equation}
\end{lemma}
\begin{proof}
    Note that $((\EB[(\bM_{\frac{k}{2^n}}-\bM_{\frac{k-1}{2^n}})^2|{\sF_{\frac{k-1}{2^n}}}]-(\qv{\bM}_{\frac{k}{2^n}}-\qv{\bM}_{\frac{k-1}{2^n}})))_{k\in\NB}$ is a martingale difference sequence. Thus, 
    \begin{equation}
    \begin{aligned}
       &\EB\norm{\sum_{k=0}^{\infty}\prn{\EB\brk{(\bM_{\frac{k}{2^n}}-\bM_{\frac{k-1}{2^n}})^2\big|{\sF_{\frac{k-1}{2^n}}}}-(\qv{\bM}_{\frac{k}{2^n}}-\qv{\bM}_{\frac{k-1}{2^n}})}}_F^2\\
       =&\sum_{k=0}^{\infty}\EB\norm{\prn{\EB\brk{(\bM_{\frac{k}{2^n}}-\bM_{\frac{k-1}{2^n}})^2\big|{\sF_{\frac{k-1}{2^n}}}}-(\qv{\bM}_{\frac{k}{2^n}}-\qv{\bM}_{\frac{k-1}{2^n}})}}_F^2\\
       \le &2\sum_{k=0}^{\infty}\EB\norm{\EB\brk{(\bM_{\frac{k}{2^n}}-\bM_{\frac{k-1}{2^n}})^2\big|{\sF_{\frac{k-1}{2^n}}}}}_F^2+2\sum_{k=0}^{\infty}\EB\norm{\qv{\bM}_{\frac{k}{2^n}}-\qv{\bM}_{\frac{k-1}{2^n}}}_F^2\\
       \le &2\sum_{k=0}^{\infty}\EB\norm{(\bM_{\frac{k}{2^n}}-\bM_{\frac{k-1}{2^n}})^2}_F^2+2\sum_{k=0}^{\infty}\EB\norm{\qv{\bM}_{\frac{k}{2^n}}-\qv{\bM}_{\frac{k-1}{2^n}}} \tr\prn{\qv{\bM}_{\frac{k}{2^n}}-\qv{\bM}_{\frac{k-1}{2^n}}}\\
       \le &2\sum_{k=0}^{\infty}\EB\tr\prn{\bM_{\frac{k}{2^n}}-\bM_{\frac{k-1}{2^n}}}^4+2\sum_{k=0}^{\infty}\EB\prn{\sup_{k\in\NB}\norm{\qv{\bM}_{\frac{k}{2^n}}-\qv{\bM}_{\frac{k-1}{2^n}}}} \tr\prn{\qv{\bM}_{\frac{k}{2^n}}-\qv{\bM}_{\frac{k-1}{2^n}}}\\
       =&2\EB\tr\sum_{k=0}^{\infty}(\bM_{\frac{k}{2^n}}-\bM_{\frac{k-1}{2^n}})^4+2\EB\prn{\sup_{k\in\NB}\norm{\qv{\bM}_{\frac{k}{2^n}}-\qv{\bM}_{\frac{k-1}{2^n}}}} \tr(\qv{\bM}_{\infty}),
    \end{aligned}
\end{equation}
where the first term converges to $0$ by Lemma~\ref{lem:fourth_moment}, and the second term converges to $0$ by the dominated convergence theorem.
\end{proof}

\begin{lemma}
\begin{equation}
    \begin{aligned}
       \prn{\norm{\sum_{k=0}^{\infty}\EB\brk{(\bM_{\frac{k}{2^n}}-\bM_{\frac{k-1}{2^n}})^2\big|{\sF_{\frac{k-1}{2^n}}}}}^{\frac{p}{2}}}_{n\in\NB}
    \end{aligned}
\end{equation}
is uniformly integrable.
\end{lemma}
\begin{proof}
Take $q>\frac{p}{2}$. It suffices to prove 
\begin{equation}
    \begin{aligned}
       \prn{\EB\norm{\sum_{k=0}^{\infty}\EB\brk{(\bM_{\frac{k}{2^n}}-\bM_{\frac{k-1}{2^n}})^2\big|{\sF_{\frac{k-1}{2^n}}}}}^{q}}_{n\in\NB}
    \end{aligned}
\end{equation}
is uniformly bounded.
\begin{equation}
    \begin{aligned}
       \EB\norm{\sum_{k=0}^{\infty}\EB\brk{(\bM_{\frac{k}{2^n}}-\bM_{\frac{k-1}{2^n}})^2\big|{\sF_{\frac{k-1}{2^n}}}}}^{q}&\le \EB\prn{\tr\prn{\sum_{k=0}^{\infty}\EB\brk{(\bM_{\frac{k}{2^n}}-\bM_{\frac{k-1}{2^n}})^2\big|{\sF_{\frac{k-1}{2^n}}}}}}^{q}\\
       &=\EB\prn{\sum_{i,j=1}^{d}\sum_{k=0}^{\infty}\EB\brk{([\bM_{\frac{k}{2^n}}]_{ij}-[\bM_{\frac{k-1}{2^n}}]_{ij})^2\big|{\sF_{\frac{k-1}{2^n}}}}}^{q}\\
       &\le d^{2(q-1)}\sum_{i,j=1}^{d}\EB\prn{\sum_{k=0}^{\infty}\EB\brk{\prn{[\bM_{\frac{k}{2^n}}]_{ij}-[\bM_{\frac{k-1}{2^n}}]_{ij}}^2\big|{\sF_{\frac{k-1}{2^n}}}}}^{q}.
    \end{aligned}
\end{equation}
Apply scalar Rosenthal-Burkholder inequality \citep{burkholder1973distribution} to $([\bM_k^{(n)}]_{ij})_{k\in\NB}$:
\begin{equation}
    \begin{aligned}
       \EB\prn{\sum_{k=0}^{\infty}\EB\brk{\prn{[\bM_{\frac{k}{2^n}}]_{ij}-[\bM_{\frac{k-1}{2^n}}]_{ij}}^2\big|{\sF_{\frac{k-1}{2^n}}}}}^{q}\le C_q K^{2q}.
    \end{aligned}
\end{equation}
Thus 
\begin{equation}
    \begin{aligned}
       \prn{\EB\norm{\sum_{k=0}^{\infty}\EB\brk{(\bM_{\frac{k}{2^n}}-\bM_{\frac{k-1}{2^n}})^2\big|{\sF_{\frac{k-1}{2^n}}}}}^{q}}_{n\in\NB}
    \end{aligned}
\end{equation}
is uniformly bounded by $d^{2q}C_q K^{2q}$.
\end{proof}
By Lemma~\ref{lem:fourth_moment} and Lemma~\ref{lem:dis_quadratic_converge}, $\lim_{n\to \infty}\sum_{k=1}^{2^n}\EB\brk{(\bM_{\frac{k}{2^n}T}-\bM_{\frac{k-1}{2^n}T})^2\big|{\sF_{\frac{k-1}{2^n}T}}}=\qv{\bM}_T$ in $L^{\frac{p}{2}}$.

Now suppose that $(\bM_t)$ is a continuous local martingale that is not necessarily bounded. We use the standard localization technique. Define a sequence of stopping times $\tau_m={\rm inf}\{t\ge 0:\norm{\bM_t}\ge m \text{ or } \norm{\qv{\bM}_t}\ge m\}$. Then $\bM_t^{(m)}:=\bM_{t\wedge\tau\wedge\tau_m}$ is a proper martingale.
By the previous results, we have
\begin{equation}
    \begin{aligned}
       \norm{\sup_{t\geq 0}\|\bM_{t\wedge\tau\wedge\tau_m}\|}_p\leq C_{R,1}\sqrt{p\vee \log d}\Big\|\norm{\qv{\bM}_{\tau\wedge\tau_m}}^{1/2} \Big\|_p \le C_{R,1}\sqrt{p\vee \log d}\Big\|\norm{\qv{\bM}_\tau}^{1/2} \Big\|_p.
    \end{aligned}
\end{equation}
Let $m\to\infty$ and we finish the proof.
\end{proof}

\subsection{Proof of Lemma~\ref{lem:truncated_phi_ineq}}\label{appendix:proof_truncated_phi_inequality}
\begin{proof}
    Consider the Lebesgue-Stieltjes integral, we have
    \begin{equation*}
        \begin{aligned}
            \EB\brk{\Phi(Y)}=&\EB\brk{\Phi(Y)\ind\brc{Y\leq a}}+\EB\brk{\Phi(Y)\ind\brc{Y> a}}\\
            \leq& \Phi(a)+\EB\brk{\Phi(Y)\ind\brc{Y> a}},
        \end{aligned}
    \end{equation*}
    and by the good $\lambda$ inequality, we further have
    \begin{equation*}
        \begin{aligned}
            \EB\brk{\Phi(Y)\ind\brc{Y> a}}=&\EB\brk{\Phi(\beta \beta^{-1}Y)\ind\brc{Y> a}}\\
            \leq&\gamma\EB\brk{\Phi( \beta^{-1}Y)\ind\brc{Y> a}}\\
            =&\gamma\int_{{a}/{\beta}}^\infty\PB\prn{Y>\beta\lambda}\Phi(d\lambda)\\
            =&\gamma\int_{{a}/{\beta}}^\infty\brk{\PB\prn{Y>\beta\lambda,Z\leq\lambda}+\PB\prn{Y>\beta\lambda,Z>\lambda}}\Phi(d\lambda)\\
            \leq&\gamma\int_{{a}/{\beta}}^\infty\brk{\epsilon\PB\prn{Y>\lambda}+\PB\prn{Z>\lambda}}\Phi(d\lambda)\\
            \leq& \gamma\epsilon\EB\brk{\Phi(Y)\ind\brc{Y> {a}/{\beta}}}+\gamma\EB\brk{\Phi(Z)}\\
            \leq& \gamma\epsilon\EB\brk{\Phi(Y)\ind\brc{Y>a}}+\gamma\epsilon\Phi(a)+\gamma\EB\brk{\Phi(Z)}\\
            \leq&\frac{\gamma\epsilon}{1-\gamma\epsilon}\Phi(a)+\frac{\gamma}{1-\gamma\epsilon}\EB\brk{\Phi(Z)}.
        \end{aligned}
    \end{equation*}
    hence
    \begin{equation*}
        \EB\brk{\Phi(Y)}\leq\frac{1}{1-\gamma\epsilon}\Phi(a)+\frac{\gamma}{1-\gamma\epsilon}\EB\brk{\Phi(Z)},
    \end{equation*}
    which is the desired conclusion.
\end{proof}

\subsection{Remaining Proof of Theorem~\ref{thm:martingale_rosenthal}}\label{appendix:omitted_proof_burkholder}
We continue from the Step~2 in Section~\ref{Subsection:proof_outline_burkholder}.
In this section, we only prove the dimension-dependent version of the inequality. 
For the dimension-free version, it can be obtained by repeating the proof line by line as presented in this section, and we omit it for brevity.

First, let us recall some important concepts. 
The first one is the \textit{tangent sequence}, which is used to construct conditionally symmetric martingales. 
The second is the \textit{conditional independence (CI) condition}; the additional random components introduced in constructing conditionally symmetric martingales will satisfy this condition as shown in Theorem~\ref{thm:constructing_CI_tangent_sequence}, and sequences satisfying the CI condition can be treated as independent sequences.
\begin{definition}\citep[Tangent Sequences, Definition~2.1]{hitczenko1990best}
    Two $(\sF_k)_{k\in\NB}$-adapted sequences $(\bX_k)_{k\in\NB}$ and $(\bY_k)_{k\in\NB}$ are said to be tangent if the conditional distributions of $\bX_k$ and $\bY_k$ given $\sF_{k-1}$ coincide almost surely for any $k\in\NB$.
\end{definition}
From this definition, it follows that $(\bM_k)_{k\in\NB}$ is conditionally symmetric if and only if $(\bX_k)_{k\in\NB}$ and $(-\bX_k)_{k\in\NB}$ are tangent.
Furthermore, suppose that there exists a sequence $(\bY_k)_{k\in\NB}$ tangent to $(\bX_k)_{k\in\NB}$, and $\bX_k, \bY_k$ are $\sF_{k-1}$-conditionally independent for any $k\in\NB$, it holds that $(\sum_{j=0}^k(\bX_j-\bY_j))_{k\in\NB}$ is a conditionally symmetric martingale.
\begin{definition}\citep[Conditional Independence (CI) Condition, Definition~2.2]{hitczenko1990best}
    An $(\sF_k)_{k\in\NB}$-adapted sequences $(\bY_k)_{k\in\NB}$ satisfies (CI) condition if there exists a $\sigma$-algebra $\sG\subset \sF$ such that the conditional distributions of $\bY_k$ given $\sF_{k-1}$ and $\sG$ coincide almost surely for any $k\in\NB$, and the sequence $(\bY_k)_{k\in\NB}$ is $\sG$-conditionally independent.
\end{definition}
The following lemma guarantees the existence of the desired sequence $(\bY_k)_{k\in\NB}$.
\begin{lemma}\citep[Lemma~2.3]{hitczenko1990best}\label{thm:constructing_CI_tangent_sequence}
    In some enlarged probability space of $(\Omega,\sF,\PB)$, there exists a sequence $(\bY_k)_{k\in\NB}$ which is tangent to $(\bX_k)_{k\in\NB}$, satisfies the (CI) condition with $\sigma$-algebra $\sG$, and satisfies $\bX_k, \bY_k$ are $\sF_{k-1}$-conditionally independent for any $k\in\NB$.
\end{lemma}
The statement and proof of this lemma require rigorous measure-theoretic language (such as regular conditional probabilities). 
To avoid introducing excessive notations, we omit the details and merely explain the idea of constructing $(\bY_k)_{k\in\NB}$ in words. 
Interested readers may refer to \citet[Lemma~2.3]{hitczenko1990best}. 
For simplicity, from now on, we shall not change the notations even though we will work in the enlarged probability space.
At time $0$, we generate an independent copy $\bY_0$ of $\bX_0$ (in an enlarged probability space to accommodate the sequence $(\bY_k)_{k\in\NB}$). 
At time $k\geq 1$, given $\sF_{k-1}$, we generate a copy $\bY_k$ that is conditionally independent of $\bX_k$. 
It is easy to see that $(\bY_k)_{k\in\NB}$ is tangent to $(\bX_k)_{k\in\NB}$ and satisfies the (CI) condition with $\sG = \sF$. 

Based on the sequence $(\bY_k)_{k\in\NB}$ given in Theorem~\ref{thm:constructing_CI_tangent_sequence}, we introduce some notations.
For any $k\in\NB$, we denote $\bN_k=\sum_{j=0}^k \bY_j$, $\tilde{\bX}_k=\bX_k-\bY_k$, $\tilde{\bM}_k=\sum_{j=0}^k \tilde{\bX}_k=\bM_k-\bN_k$. 
Then $(\bN_k)_{k\in\NB}$ is a martingale with (CI) increments, and we denote by $(\langle\bN\rangle_k)_{k\in\NB}$ the related quadratic variation process. 
In addition, $(\tilde{\bM}_k)_{k\in\NB}$ is a conditionally symmetric martingale, and we denote by $(\langle\tilde\bM\rangle_k)_{k\in\NB}$ the related quadratic variation process. 

Because $\bX_k, \bY_k$ are $\sF_{k-1}$-conditionally independent, we have $\langle\tilde\bM\rangle_k=\langle\bM\rangle_k+\langle\bN\rangle_k$ almost surely.
In addition, because $(\bY_k)_{k\in\NB}$ is tangent to $(\bX_k)_{k\in\NB}$, we have $\langle\bM\rangle_k=\langle\bN\rangle_k$ almost surely, thus $\langle\tilde\bM\rangle_k=2\langle\bM\rangle_k$ almost surely.

By Minkowski's inequality, we have
\begin{equation*}
    \norm{\sup_{k\in\NB}\|\bM_k\|}_p\leq \norm{\sup_{k\in\NB}\|\tilde\bM_k\|}_p+\norm{\sup_{k\in\NB}\|\bN_k\|}_p. 
\end{equation*}
According to Rosenthal-Burkholder inequality for conditionally symmetric martingale (Eqn.~\eqref{eq:rosenthal_burkholder_symmetric}), we have
\begin{equation*}
\begin{aligned}
    \norm{\sup_{k\in\NB}\|\tilde\bM_k\|}_p\leq& 5.01\prn{1+e^{-p/c}+\frac{1}{c\vee\log d}}\prn{\sqrt{c\vee\log d}\ e^{p/c} \norm{\norm{\qv{\tilde\bM}_\infty}^{1/2}}_p+(c\vee \log d)\norm{\sup_{k\in\NB}\norm{\tilde\bX_k}}_p },
\end{aligned}
\end{equation*}
where 
\begin{equation*}
\begin{aligned}
    \norm{\norm{\qv{\tilde\bM}_\infty}^{1/2}}_p=\sqrt{2}\norm{\norm{\qv{\bM}_\infty}^{1/2}}_p,
\end{aligned}
\end{equation*}
and
\begin{equation*}
   \norm{\sup_{k\in\NB}\norm{\tilde\bX_k}}_p \leq \norm{\sup_{k\in\NB}\norm{\bX_k}}_p+\norm{\sup_{k\in\NB}\norm{\bY_k}}_p.
\end{equation*}
To give an upper bound for $\|\sup_{k\in\NB}\|\bY_k\|\|_p$, we need the following lemma.
\begin{lemma}\citep[Lemma~1]{hitczenko1988comparison}\label{lem:bound_tangent_sequence}
   Suppose that $A_k, B_k\in\sF_k$ satisfy $\PB(A_k\mid \sF_{k-1})=\PB(B_k\mid \sF_{k-1})$ for any $k\in\NB$.
   Then it holds that
   \begin{equation*}
       \PB\prn{\bigcup_{k\in\NB}A_k}\leq 2\PB\prn{\bigcup_{k\in\NB}B_k}.
   \end{equation*}
\end{lemma}
For any $t>0$, we define the events $A_k=\{\|\bY_k\|\geq t \}$, $B_k=\{\|\bX_k\|\geq t \}$ which satisfy $\PB(A_k\mid \sF_{k-1})=\PB(B_k\mid \sF_{k-1})$ because $(\bY_k)_{k\in\NB}$ is tangent to $(\bX_k)_{k\in\NB}$.
By Lemma~\ref{lem:bound_tangent_sequence}, we obtain
\begin{equation*}
\begin{aligned}
    \PB\prn{\sup_{k\in\NB}\norm{\bY_k}\geq t}=&\PB\prn{\bigcup_{k\in\NB}A_k}\leq 2\PB\prn{\bigcup_{k\in\NB}B_k}=2\PB\prn{\sup_{k\in\NB}\norm{\bX_k}\geq t},
\end{aligned}
\end{equation*}
hence
\begin{equation*}
    \begin{aligned}
        \norm{\sup_{k\in\NB}\norm{\bY_k}}_p^p=&\EB\brk{\sup_{k\in\NB}\norm{\bY_k}^p}\\
        =&\int_0^\infty \PB\prn{\sup_{k\in\NB}\norm{\bY_k}^p\geq u} du\\
        =&p\int_0^\infty \PB\prn{\sup_{k\in\NB}\norm{\bY_k}\geq t}t^{p-1} dt\\
        \leq&2p\int_0^\infty \PB\prn{\sup_{k\in\NB}\norm{\bX_k}\geq t}t^{p-1} dt\\
        =&2\norm{\sup_{k\in\NB}\norm{\bX_k}}_p^p,
    \end{aligned}
\end{equation*}
which implies $\|\sup_{k\in\NB}\|\bY_k\|\|_p\leq2^{1/p}\|\sup_{k\in\NB}\|\bX_k\|\|_p$.

Therefore, we have
\begin{equation}\label{eq:matrix_rosenthal_burk_part_1}
\begin{aligned}
    \norm{\sup_{k\in\NB}\|\tilde\bM_k\|}_p\leq& 5.01\prn{1{+}e^{-p{/}c}{+}\frac{1}{c{\vee}\log d}}\\
    &\cdot \prn{\sqrt{2(c\vee\log d)}\ e^{p/c} \norm{\norm{\qv{\bM}_\infty}^{1/2}}_p+(1+2^{1{/}p})(c\vee \log d)\norm{\sup_{k\in\NB}\norm{\bX_k}}_p }.
\end{aligned}
\end{equation}

We only need to give an upper bound for $\norm{\sup_{k\in\NB}\|\bN_k\|}_p$.
Since $(\bY_k)_{k\in\NB}$ satisfies (CI) condition, the proof strategy is almost identical to that of an independent sum. 
The only difference is that we need to replace the expectation with the mathematical expectation conditional on $\sG$.
The key tool is the following maximal version of the symmetrization inequality \citep[Lemma~6.3]{ledoux2013probability}.
\begin{lemma}[Maximal Symmetrization Inequality]\label{lem:maximal_symm_ineq}
Let $F\colon \RB_+\to\RB_+$ be a non-decreasing convex function, and $(\gX, \|\cdot\|)$ be a separable Banach space.
Suppose that $(X_k)_{k\in\NB}$ is a sequence of independent $\gX$-valued mean zero random elements such that $\EB[F(\sup_{k\in\NB}\|\sum_{j=0}^k X_j\|)]<\infty$ for any $k\in\NB$, and $(\epsilon_k)_{k\in\NB}$ is a sequence of i.i.d. Rademacher random variables ($\PB(\epsilon_1=1)=\PB(\epsilon_1=-1)=1/2$) independent of $(X_k)_{k\in\NB}$.
Then it holds that
\begin{equation*}
    \EB\brk{F\prn{\sup_{k\in\NB}\norm{\sum_{j=0}^k X_j}}}\leq \EB\brk{F\prn{2\sup_{k\in\NB}\norm{\sum_{j=0}^k \epsilon_jX_j}}}.
\end{equation*}
\end{lemma}
The proof can be found in Appendix~\ref{appendix:proof_maximal_symmetric_ineq}.

By applying Lemma~\ref{lem:maximal_symm_ineq} with $F(x)=x^p$, we have
\begin{equation*}
    \begin{aligned}
        \norm{\sup_{k\in\NB}\|\bN_k\|}_p^p=&\EB\brk{\sup_{k\in\NB}\norm{ \sum_{j=0}^k\bY_j}^p }\\
        =&\EB\brk{\EB\brk{\sup_{k\in\NB}\norm{ \sum_{j=0}^k\bY_j}^p\Bigg| \sG }}\\
        \leq&  2^p   \EB\brk{\EB\brk{\sup_{k\in\NB}\norm{ \sum_{j=0}^k\epsilon_j\bY_j}^p\Bigg| \sG }}\\
        =&2^p\EB\brk{\sup_{k\in\NB}\norm{ \sum_{j=0}^k\epsilon_j\bY_j}^p}\\
        =&2^p\norm{\sup_{k\in\NB}\|\tilde\bN_k\|}_p^p,
    \end{aligned}
\end{equation*}
where we used the fact that $(\bY_k)_{k\in\NB}$ satisfies (CI) condition and we denoted $\tilde\bN_k=\sum_{j=0}^k\epsilon_j\bY_j$ which is also a conditionally symmetric martingale with $\langle\tilde \bN \rangle_k=\langle\bN \rangle_k=\langle\bM \rangle_k$ and $\sup_{k\in\NB}\norm{\epsilon_k\bY_k}=\sup_{k\in\NB}\norm{\bY_k}$ which implies $\|\sup_{k\in\NB}\|\epsilon_k\bY_k\|\|_p=\|\sup_{k\in\NB}\|\bY_k\|\|_p\leq2^{1/p}\|\sup_{k\in\NB}\|\bX_k\|\|_p$ for any $k\in\NB$.

Applying the Rosenthal-Burkholder inequality for conditionally symmetric martingale (Eqn.~\eqref{eq:rosenthal_burkholder_symmetric}) again, we have
\begin{equation}\label{eq:matrix_rosenthal_burk_part_2}
    \begin{aligned}
        \norm{\sup_{k\in\NB}\|\bN_k\|}_p\leq& 2\norm{\sup_{k\in\NB}\|\tilde\bN_k\|}_p\\
        \leq& 10.02\prn{1+e^{-p/c}+\frac{1}{c\vee\log d}}\prn{\sqrt{c\vee\log d}\ e^{p/c} \norm{\norm{\qv{\tilde\bN}_\infty}^{1/2}}_p+(c\vee \log d)\norm{\sup_{k\in\NB}\norm{\epsilon_k\bY_k}}_p }  \\
        \leq& 10.02\prn{1+e^{-p/c}+\frac{1}{c\vee\log d}}\prn{\sqrt{c\vee\log d}\ e^{p/c} \norm{\norm{\qv{\bM}_\infty}^{1/2}}_p+2^{1/p}(c\vee \log d)\norm{\sup_{k\in\NB}\norm{\bX_k}}_p }. 
    \end{aligned}
\end{equation}
Combining Eqn.~\eqref{eq:matrix_rosenthal_burk_part_1} with Eqn.~\eqref{eq:matrix_rosenthal_burk_part_2}, we arrive at the conclusion.
\begin{equation*}
    \begin{aligned}
        \norm{\sup_{k\in\NB}\|\bM_k\|}_p\leq& \norm{\sup_{k\in\NB}\|\tilde\bM_k\|}_p+\norm{\sup_{k\in\NB}\|\bN_k\|}_p  \\
        \leq& 5.01\prn{1{+}e^{-p{/}c}+\frac{1}{c{\vee}\log d}}\\
        &\cdot\prn{(2{+}\sqrt{2})\sqrt{c{\vee}\log d} e^{p{/}c} \norm{\norm{\qv{\bM}_\infty}^{1{/}2}}_p+(1{+}3{\cdot}2^{1{/}p})(c{\vee} \log d)\norm{\sup_{k\in\NB}\norm{\bX_k}}_p }. 
    \end{aligned}
\end{equation*}
When we take $c=p$, we determine the constants in the inequality:
The coefficient of the quadratic variation term is
\begin{equation*}
    5.01(1+e^{-1}+\frac{1}{p\vee \log d})(2+\sqrt{2})e\leq 87,
\end{equation*}
when $p\vee \log d\geq 117$, the constant can be taken as $64$.
The coefficient of the martingale difference term is
\begin{equation*}
    5.01(1+e^{-1}+\frac{1}{p\vee \log d})(1+3\cdot 2^{1/p})\leq 50,
\end{equation*}
when $p\geq 117$, the constant can be taken as $28$.

\subsection{Proof of Lemma~\ref{lem:maximal_symm_ineq}}\label{appendix:proof_maximal_symmetric_ineq}
\begin{proof}
We denote by $(X_k^\prime)_{k\in\NB}$ an independent copy of $(X_k)_{k\in\NB}$, which is also independent of $(\epsilon_k)_{k\in\NB}$. 
And we denote $\tilde{X}_k=X_k-X_k^\prime$.
Then $\epsilon_k\tilde{X}_k$ and $\tilde{X}_k$ have the same distribution.
Hence
\begin{equation*}
    \begin{aligned}
        \EB\brk{F\prn{\sup_{k\in\NB}\norm{\sum_{j=0}^k X_j}}}=&\EB\brk{F\prn{\sup_{k\in\NB}\norm{\sum_{j=0}^k \prn{X_j-\EB[X_j^\prime]}}}}\\
        \leq&\EB\brk{F\prn{\EB\brk{{\sup_{k\in\NB}\norm{\sum_{j=0}^k\tilde{X}_j}\Bigg| (X_k^\prime)_{k\in\NB}}}}}\\
        \leq& \EB\brk{F\prn{\sup_{k\in\NB}\norm{\sum_{j=0}^k \tilde{X}_j}}}\\
        =&\EB\brk{F\prn{\sup_{k\in\NB}\norm{\sum_{j=0}^k \epsilon_j\tilde{X}_j}}}\\
        \leq&\EB\brk{F\prn{\sup_{k\in\NB}\norm{\sum_{j=0}^k \epsilon_jX_j}+\sup_{k\in\NB}\norm{\sum_{j=0}^k \epsilon_jX_j^\prime}}}\\
        \leq&\frac{1}{2}\brc{\EB\brk{F\prn{2\sup_{k\in\NB}\norm{\sum_{j=0}^k \epsilon_jX_j}}}+\EB\brk{F\prn{2\sup_{k\in\NB}\norm{\sum_{j=0}^k \epsilon_jX_j^\prime}}}}\\
        =&\EB\brk{F\prn{2\sup_{k\in\NB}\norm{\sum_{j=0}^k \epsilon_jX_j}}},
    \end{aligned}
\end{equation*}
where we used the fact that the function $(x_k)_{k\in\NB}\mapsto\sup_{k\in\NB}\norm{x_k}$ from $\gX^\NB$ to $\RB_+$ is a convex function and Jensen's inequality.
\end{proof}

%% file: proof_markov.tex
In this section, we complete the proofs omitted in Section~\ref{Subsection:proof_outline_markov}.

\subsection{Proof of Lemma~\ref{lem:poisson_equation}}\label{appendix:proof_poisson_equation}
\begin{proof}
To derive the desired conclusion, we need the following useful lemma, whose proof is in Appendix~\ref{appendix:proof_integral_bound_tv}.
\begin{lemma}\label{lem:integral_bound_tv}
For any $\xi_1, \xi_2\in\gP$, vector-valued function $\bh\colon\gZ\to\CB^{d}$, and matrix-valued function $\bH\colon\gZ\to\CB^{d_1{\times}d_2}$, it holds that
  \begin{enumerate}
    \item Assume $\norm{\norm{\bh}}_{\infty}, \norm{\norm{\bH}}_{\infty}<\infty$. Then it holds that
    \begin{equation*}
    \begin{aligned}
        &\norm{\xi_1(\bh)-\xi_2(\bh)}\leq \norm{\norm{\bh}}_{\infty}\norm{\xi_1-\xi_2}_{\TV},\\
        &\norm{\xi_1(\bH)-\xi_2(\bH)}\leq \norm{\norm{\bH}}_{\infty}\norm{\xi_1-\xi_2}_{\TV}.
    \end{aligned}
    \end{equation*}
    \item Assume $\norm{\norm{\bh}}_{V^\alpha}, \norm{\norm{\bH}}_{V^\alpha}<\infty$ for some $\alpha\in(0,1]$. Then it holds that
    \begin{equation*}
    \begin{aligned}
        &\norm{\xi_1(\bh)-\xi_2(\bh)}\leq \norm{\norm{\bh}}_{V^\alpha}\norm{\xi_1-\xi_2}_{V^\alpha},\\
        &\norm{\xi_1(\bH)-\xi_2(\bH)}\leq \norm{\norm{\bH}}_{V^\alpha}\norm{\xi_1-\xi_2}_{V^\alpha}.
    \end{aligned}
    \end{equation*}
\end{enumerate}  
\end{lemma}

Equipped with Lemma~\ref{lem:integral_bound_tv}, we are ready to prove this lemma.
\begin{enumerate}[(a)]
    \item  
    By \citet[Proposition~21.2.3]{douc2018markov}, we only need to show that $\|\bG(z)\|<\infty$ for any $z\in\gZ$.
    \begin{equation*}
        \begin{aligned}
            \norm{\norm{\bG}}_{\infty}=&\sup_{z\in\gZ}\norm{\sum_{k=0}^\infty\prn{\rQ^k \bF}(z)}\\
            \leq&\sup_{z\in\gZ}\sum_{k=0}^\infty\norm{\prn{\rQ^k \bF}(z)}\\
            =&\sup_{z\in\gZ}\sum_{k=0}^\infty\norm{\rQ^k(z,\cdot)(\bF)-\mu(\bF) }\\
            \leq&\sup_{z\in\gZ}\sum_{k=0}^\infty\norm{\norm{\bF}}_{\infty}\norm{\rQ^k(z,\cdot)-\mu}_{\TV}\\
            \leq&\norm{\norm{\bF}}_{\infty}\sum_{k=0}^\infty\sup_{z\in\gZ}\norm{\rQ^k(z,\cdot)-\mu}_{\TV}\\
            \leq&2\norm{\norm{\bF}}_{\infty}\sum_{k=0}^\infty\prn{\frac{1}{4}}^{\lfloor k/\tmix\rfloor}\\
            \leq&2\tmix \norm{\norm{\bF}}_{\infty}\sum_{k=0}^\infty\prn{\frac{1}{4}}^{k}\\
            =&\frac{8}{3}\tmix \norm{\norm{\bF}}_{\infty}.
        \end{aligned}
    \end{equation*}
    where we used Lemma~\ref{lem:integral_bound_tv} and Assumption~\ref{assumption:uniform_geo_ego}.

    If we further assume $\bF^2(z)\preccurlyeq\bUp$ for all $z\in\gZ$, for some $\bUp\in\HB^d_+$, then we have the following result: for any $\bu\in\CB^d$ and $z\in\gZ$, it holds that
    \begin{equation*}
        \begin{aligned}
            \norm{\bG(z)\bu}=&\norm{\sum_{k=0}^\infty\prn{\rQ^k \bF}(z)\bu}\\
            \leq&\sum_{k=0}^\infty\norm{\prn{\rQ^k \bF}(z)\bu}\\
            =&\sum_{k=0}^\infty\norm{\rQ^k(z,\cdot)(\bF\bu)-\mu(\bF\bu) }\\
            \leq&\sum_{k=0}^\infty \norm{\norm{\bF\bu}}_{\infty}\norm{\rQ^k(z,\cdot)-\mu}_{\TV}\\
            \leq&\frac{8}{3}\tmix \norm{\norm{\bF\bu}}_{\infty},
        \end{aligned}
    \end{equation*}
    where we used Lemma~\ref{lem:integral_bound_tv}.
    Hence
    \begin{equation*}
        \begin{aligned}
            \bu^* \bG^2(z)\bu =&\norm{\bG(z)\bu}^2\\
            \leq& \frac{64}{9}\tmix \norm{\norm{\bF\bu}}_{\infty}^2\\
            =&\frac{64}{9}\tmix\sup_{z^\prime\in\gZ}\bu^* \bF^2(z^\prime)\bu\\
            \leq&\frac{64}{9}\tmix\bu^* \bUp\bu,
        \end{aligned}
    \end{equation*}
    which is the desired conclusion.

    \item 
    We first state a lemma, which is a immediate corollary of Jensen's inequality.
    \begin{lemma}\citep[Lemma~12]{durmus2023rosenthal}\label{lem:V_alpha_ergo}
    Assume Assumption~\ref{assumption:V_ego}. Then for any $\alpha\in(0,1]$, $k\in\NB$, and $\xi\in\gP_V$, it holds that
    \begin{equation*}
        \norm{\xi\rQ^k-\mu}_{V^\alpha}\leq 2^{1-\alpha}\varkappa^\alpha\prn{\mu(V)+\xi(V)}^{\alpha}(1/4)^{{\lfloor k/\tmix\rfloor}\alpha}.
    \end{equation*}
    \end{lemma}
        
    The proof is almost identical:    
    \begin{equation*}
        \begin{aligned}
            \norm{\bG(z)}=&\norm{\sum_{k=0}^\infty\prn{\rQ^k \bF}(z)}\\
            \leq&\sum_{k=0}^\infty\norm{\rQ^k(z,\cdot)(\bF)-\mu(\bF) }\\
            \leq&\sum_{k=0}^\infty\norm{\norm{\bF}}_{V^\alpha}\norm{\rQ^k(z,\cdot)-\mu}_{V^\alpha}\\
            \leq&2^{1-\alpha}\varkappa^\alpha\prn{\mu(V)+V(z)}^{\alpha}\norm{\norm{\bF}}_{V^\alpha}\sum_{k=0}^\infty(1/4)^{{\lfloor k/\tmix\rfloor}\alpha}\\
            \leq&\brk{\varkappa^\alpha2^{3-\alpha}/(3\alpha)}(\mu(V)+V(z))^{\alpha}\tmix\norm{\|\bF\|}_{V^\alpha},
        \end{aligned}
    \end{equation*}
    where we used Lemma~\ref{lem:integral_bound_tv}, Lemma~\ref{lem:V_alpha_ergo}, and the estimate in \citet[Eqn~(146)]{durmus2023rosenthal} in the last inequality.
\end{enumerate}
\end{proof}

\subsection{Proof of Lemma~\ref{lem:crude_rosenthal_markov}}\label{appendix:proof_crude_rosenthal_markov}
\begin{proof}
\begin{enumerate}[(a)]
    \item By Minkowski inequality, we have
\begin{equation*}
    \norm{\norm{\bS_n}}_{\mu,p}\le\norm{\norm{\sum_{k=0}^{\lfloor{n}/{\tmix}\rfloor\tmix-1}{\bF}(Z_k)}}_{\mu,p}+\tmix\norm{\norm{\bF}}_{\infty}.
\end{equation*}
We define the function $\bG_{\tmix}\colon \gZ\to\HB^d$ as $\bG_{\tmix}=\sum_{k=0}^\infty \rQ^{k\tmix} \bF $, 
which satisfies the Poisson equation $\bG_{\tmix}-\rQ^{\tmix} \bG_{\tmix}=\bF$ .
By Lemma~\ref{lem:poisson_equation} (a), $\bG_{\tmix}$ satisfies the Poisson equation $\bG_{\tmix}-\rQ^{\tmix} \bG_{\tmix}=\bF$ and $\norm{\norm{\bG_{\tmix}} }_{\infty}\leq  ({8}/{3})\norm{\norm{\bF} }_{\infty}$.
By Minkowski inequality, we have
\begin{equation*}
    \begin{aligned}
    \norm{\norm{\sum_{k=0}^{\lfloor\frac{n}{\tmix}\rfloor\tmix-1}{\bF}(Z_k)}}_{\mu,p}
    =&\norm{\norm{\sum_{r=0}^{\tmix-1}\sum_{k=0}^{\lfloor\frac{n}{\tmix}\rfloor-1}\prn{{\bG_{\tmix}}(Z_{k\tmix+r})-\prn{\rQ^{\tmix} \bG_{\tmix}}(Z_{k\tmix+r})}}}_{\mu,p}\\
    \le&\sum_{r=0}^{\tmix-1}\norm{\norm{\sum_{k=0}^{\lfloor\frac{n}{\tmix}\rfloor-1}\prn{{\bG_{\tmix}}(Z_{k\tmix+r})-\prn{\rQ^{\tmix} \bG_{\tmix}}(Z_{k\tmix+r})}}}_{\mu,p}\\
    \le&\sum_{r=0}^{\tmix-1}\norm{\norm{\sum_{k=1}^{\lfloor\frac{n}{\tmix}\rfloor}\prn{{\bG_{\tmix}}(Z_{k\tmix+r})-\prn{\rQ^{\tmix} \bG_{\tmix}}(Z_{(k-1)\tmix+r})}}}_{\mu,p}\\
    &+2\tmix\norm{\norm{\bG_{\tmix}} }_{\infty}.
    \end{aligned}
\end{equation*}
The summand of the first term is the spectral norm of a matrix martingale $\bM^{(\tmix)}_{\lfloor{n}/{\tmix}\rfloor}$, by Rosenthal-Burkholder inequality for matrix martingales (Theorem~\ref{thm:martingale_rosenthal}), we have
\begin{equation*}
    \begin{aligned}
    &\norm{\norm{\sum_{k=1}^{\lfloor\frac{n}{\tmix}\rfloor}\prn{{\bG_{\tmix}}(Z_{k\tmix+r})-\prn{\rQ^{\tmix} \bG_{\tmix}}(Z_{(k-1)\tmix+r})}}}_{\mu,p}\\
    \le& C_{1}\sqrt{p\vee \log d}\norm{\norm{\sum_{k=1}^{\lfloor\frac{n}{\tmix}\rfloor}\EB\brk{\prn{{\bG_{\tmix}}(Z_{k\tmix+r})-\prn{\rQ^{\tmix} \bG_{\tmix}}(Z_{(k-1)\tmix+r})}^2\big|{\sF_{r,k-1}}}}^{1/2}}_{\mu,p} \\
    &+ C_{2}(p\vee\log d)\norm{ \sup_{k\le \lfloor\frac{n}{\tmix}\rfloor}\|{\bG_{\tmix}}(Z_{k\tmix+r})-\prn{\rQ^{\tmix} \bG_{\tmix}}(Z_{(k-1)\tmix+r})\|}_{\mu,p}\\
    \le& C_{1}\sqrt{p\vee \log d}\Big\|\prn{\sum_{k=1}^{\lfloor\frac{n}{\tmix}\rfloor}\prn{2\norm{\norm{\bG_{\tmix}} }_{\infty}}^2}^{1/2} \Big\|_{\mu,p} + C_{2}(p\vee\log d)\norm{ 2\norm{\norm{\bG_{\tmix}} }_{\infty}}_{\mu,p}\\
    =&2C_{1}\sqrt{(p\vee \log d)\lfloor\frac{n}{\tmix}\rfloor}\norm{\norm{\bG_{\tmix}} }_{\infty}  + 2C_{2}(p\vee\log d) \norm{\norm{\bG_{\tmix}} }_{\infty}.
    \end{aligned}
\end{equation*}
where $\sF_{r,k-1}=\sigma(Z_{j\tmix+r},0\le j\le k-1)$. Combining the upper bounds above, we have
\begin{equation*}
    \begin{aligned}
    \norm{\norm{\bS_n}}_{\mu,p}&\leq\tmix(2C_{1}\sqrt{(p{\vee} \log d)\lfloor{n}/{\tmix}\rfloor}\norm{\norm{\bG_{\tmix}} }_{\infty}  {+} 2C_{2}(p{\vee}\log d) \norm{\norm{\bG_{\tmix}} }_{\infty}\\
    &{+}2\tmix\norm{\norm{\bG_{\tmix}} }_{\infty}{+}\tmix\norm{\norm{\bF}}_{\infty}\\
    &\le \frac{16}{3}C_{1}\sqrt{(p\vee \log d)n\tmix}\norm{\norm{\bF}}_{\infty}  + \prn{\frac{16}{3}C_{2}(p\vee\log d)+\frac{19}{3}}\tmix \norm{\norm{\bF}}_{\infty}.
    \end{aligned}
\end{equation*}
Furthermore, if there exists an $\bUp\in\HB^d_+$ such that $\bF^2(z)\preccurlyeq \bUp$ for all $z\in\gZ$, then by Lemma~\ref{lem:poisson_equation} (a), it holds that $\bG^2_{\tmix}(z)\preccurlyeq (64/9)\bUp$ for all $z\in\gZ$.
By Theorem~\ref{thm:martingale_rosenthal} (b), the above inequality still holds if we replace $d$ with $e\cdot r([(64/9)\bUp]\wedge [(1/\lfloor{n}/{\tmix}\rfloor)\|\|\langle{\bM^{(\tmix)}}\rangle_{\lfloor{n}/{\tmix}\rfloor}\|\|_{\mu,\frac{p}{2}}])$, and double the right-hand side.
    \item The proof is identical to that of part (a), except that we need to use Lemma~\ref{lem:poisson_equation} (b) instead.
\end{enumerate}

\end{proof}

\subsection{Proof of Lemma~\ref{lem:integral_bound_tv}}\label{appendix:proof_integral_bound_tv}
\begin{proof}
\begin{enumerate}
    \item 
    We first show the inequality for the vector-valued function.
    Recall the dual representation of the total variation \citep[Proposition~D.2.4]{douc2018markov}:
    \begin{equation*}
        \norm{\xi_1-\xi_2}_{\TV}=\sup_{f\colon\|f\|_\infty\leq 1}(\xi_1(f)-\xi_2(f)),
    \end{equation*}
    then by the dual representation of the Euclidean norm, we have
    \begin{equation*}
        \begin{aligned}
            \norm{\xi_1(\bh)-\xi_2(\bh)}=&\sup_{\bu\colon\norm{\bu}\leq 1}\bu^*\prn{\xi_1(\bh)-\xi_2(\bh)}\\
            =&\sup_{\bu\colon\norm{\bu}\leq 1}\prn{\xi_1(\bu^*\bh)-\xi_2(\bu^*\bh)}\\
            \leq&\sup_{\bu\colon\norm{\bu}\leq 1} \|\bu^*\bh \|_{\infty}\norm{\xi_1-\xi_2}_{\TV}\\
            \leq&\sup_{\bu\colon\norm{\bu}\leq 1} \|\bu\|\norm{\|\bh \|}_{\infty}\norm{\xi_1-\xi_2}_{\TV}\\
            =&\norm{\|\bh \|}_{\infty}\norm{\xi_1-\xi_2}_{\TV},
        \end{aligned}
    \end{equation*}
    which is desired.

    Now, we deal with the matrix case:
    \begin{equation*}
        \begin{aligned}
            \norm{\xi_1(\bH)-\xi_2(\bH)}=&\sup_{\bu\colon\norm{\bu}\leq 1}\norm{\prn{\xi_1(\bH)-\xi_2(\bH)}\bu}\\
            =&\sup_{\bu\colon\norm{\bu}\leq 1}\norm{\xi_1(\bH\bu)-\xi_2(\bH\bu)}\\
            \leq&\sup_{\bu\colon\norm{\bu}\leq 1}\norm{\|\bH\bu \|}_{\infty}\norm{\xi_1-\xi_2}_{\TV}\\
            \leq&\sup_{\bu\colon\norm{\bu}\leq 1}\norm{\bu}\norm{\|\bH \|}_{\infty}\norm{\xi_1-\xi_2}_{\TV}\\
            =&\norm{\|\bH \|}_{\infty}\norm{\xi_1-\xi_2}_{\TV},
        \end{aligned}
    \end{equation*}
    where we used the inequality for the vector-valued function.

    \item The proof is identical to (a) if we note that
    \begin{equation*}
    \begin{aligned}
        \xi_1(f)-\xi_2(f)=&\int_\gZ f(z) (\xi_1-\xi_2)(dz)\\
        =&\int_\gZ \frac{f(z)}{V^\alpha(z)}V^{\alpha}(z) (\xi_1-\xi_2)(dz)\\
        \leq&\int_\gZ \frac{\abs{f(z)}}{V^\alpha(z)}V^{\alpha}(z) \abs{\xi_1-\xi_2}(dz)\\
        \leq& \sup_{z\in\gZ} \frac{\abs{f(z)}}{V^\alpha(z)}\int_\gZ V^{\alpha}(z) \abs{\xi_1-\xi_2}(dz)\\
        =&\norm{f}_{V^\alpha}\norm{\xi_1-\xi_2}_{V^\alpha}.
    \end{aligned}
    \end{equation*}
    Then we only need to replace $\|\cdot\|_{\infty}$ by $\norm{\cdot}_{V^\alpha}$, and $\|\cdot\|_{\TV}$ by $\norm{\cdot}_{V^\alpha}$ in (a).
\end{enumerate}

\end{proof}